\documentclass{svjour3}
\journalname{}
\usepackage{amsmath,amsfonts,amssymb}
\usepackage{graphicx}
\graphicspath{{./}{./figs/}}
\usepackage[english]{babel}
\usepackage[utf8]{inputenc}

\begin{document}
\title{Implicit-explicit schemes for compressible Cahn-Hilliard-Navier-Stokes
  equations}
\author{Pep Mulet}
\institute{
  Department of Mathematics,  Universitat de
  València    (Spain); email: mulet@uv.es.
}
\maketitle

\begin{abstract}
  
  The  isentropic compressible Cahn-Hilliard-Navier-Stokes equations
  is a system of fourth-order 
  partial differential equations that model the  evolution of some
  binary fluids   under convection.

  The purpose of this paper is the design of efficient numerical schemes to
  approximate the solution of initial-boundary value problems with these
  equations. The efficiency stems from the implicit treatment  of the
  high-order terms in the equations. Our proposal is a second-order linearly
  implicit-explicit time stepping scheme applied in a method of lines
  approach, in which the convective terms are treated explicitly and only
  linear systems have to be solved.

  Some experiments are performed to assess the validity and efficiency
  of this proposal.

\end{abstract}

\section{Introduction}
According to Kynch's theory  (see \cite{Kynch52}) for sedimentation of  homogeneous
monodisperse suspensions,  consisting of  solid
spherical particles of the same diameter and density immersed in a
viscous fluid,  two interfaces form in the settling process:  a descending interface between the clear liquid
and the initial homogeneous mixture and an ascending interface between
the maximally concentrated mixture and the initial homogeneous
mixture.

In \cite{Siano79}, it is observed that, after several days, a colloidal monodisperse
suspension of polystyrene particles sediments
forming a layered structure a fact that contradicts Kynch's theory.
   A spinodal decomposition, governed by the Cahn-Hilliard
equation \cite{CahnHilliard59}, is then conjectured as the underlying
mechanism that explains this phenomenon.

The Cahn-Hilliard equation cannot explain, by itself, this
layering phenomenon, for it does not take into account the gravitational
force. This may be  introduced into the model by means of
conservation of individual species and bulk momenta. Ignoring
temperature, this yields a
system of equations, the  isentropic
Cahn-Hilliard-Navier-Stokes equations \cite{LT98,AbelsFeireisl08}, which  are a system of fourth-order
partial differential equations that model the  evolution of
mixtures of  binary fluids  under convective effects.

Although (quasi) incompressible versions of these equations might be more
suitable for explaining the cited layering phenomenon, we consider the
compressible case for the evolution of, e.g. foams, solidification
processes, fluid–gas interface. 

This paper aims to propose numerical methods that use
implicit-explicit time-stepping schemes to avoid the severe 
restriction posed by the fourth-order terms for the efficient numerical
solution of boundary-initial problems with these equations. 

The paper is organized as follows: in section \ref{s:nsch} the
compressible isentropic Cahn-Hilliard-Navier-Stokes equations are
introduced;  in section
\ref{s:2d} implicit-explicit Runge-Kutta numerical schemes for the
two-dimensional equations are proposed;  in   section \ref{s:numerics}
we perform some numerical experiments to assess the efficiency
of our proposals; finally, in section \ref{s:conclusion} we draw some
conclusions and give some perspectives for future research.

\section{Cahn-Hilliard-Navier-Stokes equations}\label{s:nsch}
The models in this exposition are based on \cite{AbelsFeireisl08}.
We denote by $c_i$,  the mass concentration of species $i=1,2$, 
by $c=c_1-c_2$, by  $\rho$  the density of the mixture and by $\boldsymbol{v} $ its
bulk velocity (we use boldface for vector variables).

We denote by $\Omega$ the open set in $\mathbb{R}^3$ that is filled by
the fluids and by $\varepsilon$ a parameter related to the thickness of the diffuse
interface of the fluid mixture. The Ginzburg-Landau  free energy in some
region $V\subseteq \Omega$ of the
immiscible compressible two-phase fluid is
\begin{align*}
  E(\rho, c)&=\int_{V}(\rho f(\rho, c)+
    \frac{\varepsilon}{2}|\nabla c|^2)dx\\
  f(\rho, c)&=f_e(\rho)+ \psi(c)
\end{align*}    
where $\psi(c)  = \frac{1}{4}(c^2-1)^2$ is a double-well
potential function and $f_{e}$ is the specific Helmholtz free energy
of an equivalent one-phase fluid.

The isentropic compressible Cahn-Hilliard-Navier-Stokes equations with
gravitation are the following equations:
\begin{equation}\label{eq:nsch}
  \left\{
    \begin{aligned}
      \rho_t+\mathop{\text{div}} \left(\rho \boldsymbol{v}\right)&=0,\\
      (\rho\boldsymbol{v})_t+\mathop{\text{div}} \left(\rho \boldsymbol{v}\otimes
        \boldsymbol{v} \right)&=\rho \boldsymbol{G}+\mathop{\text{div}} \mathbb{T}, 
      \\
      (\rho c)_t+\mathop{\text{div}} \left(\rho c\boldsymbol{v}\right)&=\Delta \mu,
    \end{aligned}
  \right.
\end{equation}
where       $\mathop{\text{div}}$ is the divergence
operator with respect to $x\in{\mathbb{R}}^3$, 
the first equation is the continuity equation for the mixture,
the second is the equation for conservation of bulk momenta. In these equations 
\begin{align*}
  &\mathbb{T}=\nu(c)(\nabla \boldsymbol{v}+\nabla
    \boldsymbol{v}^T)+(\lambda(c)\mathop{\text{div}} \boldsymbol{v} -         p(\rho, c))\mathbb{I}+\frac{\varepsilon}{2}|\nabla c|^2
    \mathbb{I}-\varepsilon ( \nabla c\otimes \nabla  c)
\end{align*}    
is the stress tensor,
$p(\rho, c)=\rho^2\frac{\partial f(\rho, c)}{\partial \rho}$ is the
fluid pressure, $\nu(c),
\lambda(c)>0$ are   the viscosity coefficients,  $\boldsymbol{G}$ is
the gravitational  acceleration,
and
\begin{align*}
  &\mu=\psi'(c) - \frac{\varepsilon}{\rho}
    \Delta c,
\end{align*}
is the chemical potential.

The $(\rho, \rho \boldsymbol{v})$-subsystem, with $\varepsilon=0$, form the compressible isentropic Navier-Stokes equations. The equation for $\rho c$, for constant $\rho$ (which may be assumed to be 1) and $\boldsymbol{v}=0$, is the Cahn-Hilliard equation \cite{CahnHilliard59}.
\begin{align}\label{eq:ch}
      c_t=\Delta \big(\psi'(c) - \varepsilon    \Delta c\big).  
\end{align}

These equations are supplemented by initial conditions $\rho_0, \boldsymbol{v}_0,
c_0$ and the boundary conditions
\begin{align}\label{eq:bc}
  \boldsymbol{v}|_{\partial \Omega}=
  \nabla c\cdot \boldsymbol{n}|_{\partial \Omega}=
  \nabla \mu\cdot \boldsymbol{n}|_{\partial \Omega}=0,
\end{align}
where $\boldsymbol{n}$ is the outward normal vector to the boundary.

In \cite{AbelsFeireisl08} it is proved that these
equations  admit  weak solutions, with
renormalization of $\rho$ in the sense of Di Perna and Lions, in any
interval $[0, T], T>0$, provided $\gamma > \frac{3}{2}$, $0\geq
\rho_0\in L^{\gamma}(\Omega), \rho_0|v_0|^2 \in L^1(\Omega), c_0\in H^1(\Omega)$.

We henceforth consider $\nu(c), \lambda(c)$ constant  and  $p=p(\rho)=\rho^{\gamma}$, for the adiabatic constant $\gamma > 1.5$, which
corresponds to
$f_e(\rho)=\frac{\rho^{\gamma-1}}{\gamma-1}$. Therefore, the equation
for the conservation of bulk momenta can be rewritten as:
\begin{multline*}
  (\rho\boldsymbol{v})_t+\mathop{\text{div}} \left(\rho \boldsymbol{v}\otimes
    \boldsymbol{v} +p(\rho)\mathbb{I}\right)
  \\
  =\rho \boldsymbol{G}+(\nu+\lambda)\nabla \mathop{\text{div}} \boldsymbol{v}+\nu \Delta    \boldsymbol{v}
    +\frac{\varepsilon}{2}\nabla |\nabla c|^2 -\varepsilon \mathop{\text{div}}( \nabla c\otimes \nabla  c).
  \end{multline*}

As expected, the $\rho$ and $q=\rho c$ variables are conserved since
the respective associated fluxes, 
\begin{align*}
  &\rho\boldsymbol{v}\cdot\boldsymbol{n},\quad  (q\boldsymbol{v}-\nabla\mu)\cdot\boldsymbol{n},\\
\end{align*}  
vanish at the boundary due to \eqref{eq:bc}.

The two-dimensional version of these equations, for $\boldsymbol{v}=(v_1, v_2)$, is:
\begin{equation}\label{eq:chns2d}
  \begin{aligned}
    \rho_t +{(\rho v_1)_x+(\rho v_2)_{y}}&=0,\\                  
    (\rho v_1)_t+{(\rho v_1^2 +\rho^\gamma )_{x}+(\rho v_1
      v_2)_{y}}&=\frac{\varepsilon}{2}(c_{y}^2-c_{x}^2)_x
    -\varepsilon
    (c_{x}c_y)_{y}
    \\
    &\qquad +{\nu \Delta v_1+(\nu+\lambda)(  (v_1)_{xx}+(v_2)_{xy})},\\
    (\rho v_2)_t+{(\rho v_1v_2)_{x}+(\rho
      v_2^2+\rho^{\gamma})_{y}}&={\rho
      G}+\frac{\varepsilon}{2}(c_{x}^2-c_{y}^2)_y
    -\varepsilon(c_xc_y)_x \\
    &\qquad+{\nu \Delta v_2  +(\nu+\lambda)(  (v_1)_{xy}+(v_2)_{yy})},\\  
    (\rho c)_t+{(\rho c v_1)_{x}+(\rho c v_2)_{y}}&={\Delta(\psi'(c) - \frac{\varepsilon}{\rho}  \Delta c)},
  \end{aligned}
\end{equation}
where $\Delta w=w_{xx}+w_{yy}$ and gravity acts along the $y$ coordinate.

We also consider the   one-dimensional version of these equations:
\begin{equation}\label{eq:1}
  \begin{aligned}
    \rho_t+ (\rho v)_x&=0,\\
    (\rho v)_t+ (\rho v^2 +p(\rho) )_x&=\rho G+
    \Big((2\nu+\lambda)v_x-\frac{\varepsilon}{2}c_x^2\Big)_x,
    \\
    (\rho c)_t+ (\rho cv)_x&=(\psi'(c) - \frac{\varepsilon}{\rho}  c_{xx})_{xx},
  \end{aligned}
\end{equation}
where gravity acts along the $x$ coordinate.

\subsection{Spinodal decomposition}\label{ss:spinodal}
The Cahn-Hilliard equation was proposed in \cite{CahnHilliard59}  (see also \cite{Elliott89}) to model the separation of a homogeneous mixture of two incompressible  fluids, the first  of them stable with respect to the presence of small quantities of the second one, and this one unstable with respect to the presence of small quantities of the first one. The boundary of the unstable region in $(c, T, p)$-space, $T, p$ being temperature and pressure, respectively, is given by the equation $\frac{\partial^2 G}{\partial c^2}=0$, where $G(c,T,p)$ is the Gibbs free energy density of the fluid, and is usually named the spinodal.

To analyze the spinodal decomposition, we consider the  linearization
of \eqref{eq:ch} about a constant state $c_0$ in the spinodal region
$(-\frac{1}{\sqrt3}, \frac{1}{\sqrt3})$,
i.e., $\psi''(c_0)<0$:
$$c(x,t)=c_0+u(x, t),$$
for assumedly small $u$, with $\int_{\Omega} u=0$. Notice then that
$\psi'(c_0+u(x, t))= \psi'(c_0)+\psi''(c_0)u(x, t)+\mathcal{O}(u^2)$,
therefore  the  linearized Cahn-Hilliard equation:
 \begin{align}\label{eq:2222}
   &u_t= \psi''(c_0)\Delta u - \varepsilon \Delta^2 u.
 \end{align}
 is deduced from the Cahn-Hilliard equation \eqref{eq:ch}.
 
 By separation of variables,  functions of the form
 \begin{align}\label{eq:solspin0}
   u(x, t)=v(t) \prod_{i=1}^{3}\cos(k_i\pi x_i),\quad
   k_i\in\mathbb{N},
 \end{align}
 satisfy the
homogeneous Neumann boundary conditions and $\int_{\Omega} u(x, t)dx =0$, and are, therefore, solutions of \eqref{eq:2222}  if 
 \begin{align*}
   &v'(t)=\Big(-\psi''(c_0) 
     \sum_{i}k_i^2\pi^2-\varepsilon
     (\sum_{i}k_i^2\pi^2)^2) v(t),
 \end{align*}
 which yield solutions
 \begin{align}\label{eq:solspin}
   &v(t)=v(0)e^{-\big(\psi''(c_0) 
     \sum_{i}k_i^2\pi^2+\varepsilon
     (\sum_{i}k_i^2\pi^2)^2\big)t}.
 \end{align}

 The linearized equation  \eqref{eq:ch} will therefore develop instabilities
 provided 
\begin{align*}
\psi''(c_0) 
     \sum_{i}k_i^2\pi^2+\varepsilon
     (\sum_{i}k_i^2\pi^2)^2< 0,
\end{align*}
for some $ k_i \geq 1$, and 
this will be so if $\varepsilon \pi^2+\psi''(c_0) < 0$. But these
instabilities, triggered by the linearized equation, will grow until
some point, when the nonlinear character of the Cahn-Hilliard
equations makes the linearization no longer valid, and will stop, since   
energy decreases with time as shown here:
\begin{align*}
  &\frac{d}{dt} F(c(\cdot, t))=
\int_{\Omega} (\psi'(c)c_t+\varepsilon\nabla c \nabla c_t) dx 
  =\int_{\Omega} (\psi'(c)c_t-\varepsilon\Delta c  c_t) dx +
  \int_{\partial\Omega} \varepsilon  c_t \nabla c \cdot \boldsymbol{n} dx \\
&  =-\int_{\Omega} (\psi'(c)-\varepsilon\Delta c  )^2 dx \leq 0,
\end{align*}
by  using differentiation under the integral, integration by parts and  the Neumann boundary conditions for $c$.

\section{Numerical schemes}\label{s:2d}
Numerical schemes for the Cahn-Hilliard equation can be found, e.g.,
in \cite{ElliottFrench87,Elliott89}, and,  for the
quasi-incompressible Cahn-Hilliard-Navier-Stokes, in, e.g.,
\cite{Jacmin99,Yue04,ShenYang10}. In \cite{HeShi20} there is a
numerical study for compressible Cahn-Hilliard-Navier-Stokes that
mainly focus on the convective part.

Our purpose is to design finite differences numerical methods for the efficient approximate solution of the two-dimensional equations in Section  \ref{s:nsch}. 
For this, we consider $\Omega=(0,1)^2$ and the equispaced computational grid given
by the $M^2$ nodes $x_{i,j}=((i-\frac{1}{2})h, (j-\frac{1}{2})h),
i,j=1,\dots,M$, where $h=\frac{1}{M}$ and  denote by $(x, y)$ the spatial variable.

We denote by 
\begin{equation*}
  u=(\rho, \boldsymbol{m}, q), \boldsymbol{m}=(m_1,m_2)=(\rho v_1,\rho v_2)=\rho \boldsymbol{v}, q=\rho c,
\end{equation*}
the vector of conserved variables and aim to approximate \eqref{eq:chns2d}
by a spatial semi-discretization consisting of $4M^2$  ordinary
differential equations
\begin{align*}
  u_{k,i,j}'(t)={\mathcal{L}}(U(t))_{k,i,j}, k=1,\dots,4,i,j=1,\dots,M,
\end{align*}
for $4M^2$ unknowns $u_{k,i,j}(t)\in{\mathbb{R}}^4$ which are approximations of
$ u_k(x_{i,j}, t)$ and form the $4M^2$ (column) vector function $U(t)$ by using
lexicographical order so that
\begin{align*}
  &U=\begin{bmatrix}
    \varrho\\
    \varrho*V_1\\
    \varrho*V_2\\
    \varrho*C
  \end{bmatrix},\quad (\varrho*S)_{i}=\varrho_iS_i,\\
  &\rho(x_{i,j},t)\approx \varrho_{M(i-1)+j}(t),\\
  &v_k(x_{i,j},t)\approx
    (V_k)_{M(i-1)+j}(t), k=1,2,  \\
  &c(x_{i,j},t)\approx C_{M(i-1)+j}(t).
\end{align*}
For the sake of notation and simplicity we
seamlessly use a slight abuse of notation when identifying, e.g.,
$ \varrho_{i,j}\equiv \varrho_{M(i-1)+j}$. We also use superindices
for $M^2$ block notation, e.g., $U^1=\varrho$.

The nonzero terms in the spatial semidiscretization
\begin{equation}\label{eq:2}
  \begin{aligned}
    \mathcal{L}(U)&={\mathcal{C}}(U)+{\mathcal{L}}_1(U)+{\mathcal{L}}_2(U)+{\mathcal{L}}_3(U)+{\mathcal{L}}_4(U)\\
  \end{aligned}    
\end{equation}
are the following:
\begin{align*}
  &{\mathcal{C}}(U)_{1,i,j}\approx -(      (\rho v_1)_x+(\rho v_2)_{y})(x_{i,j}, t),\\
  &{\mathcal{C}}(U)_{2,i,j}\approx -(       (\rho v_1^2 +\rho^\gamma )_{x}+(\rho v_1 v_2)_{y})(x_{i,j}, t),\\              
  &{\mathcal{C}}(U)_{3,i,j}\approx -(       (\rho v_1v_2)_{x}+(\rho v_2^2+\rho^{\gamma})_{y})(x_{i,j}, t),\\  
  &{\mathcal{C}}(U)_{4,i,j}\approx -(       (\rho c v_1)_{x}+(\rho c v_2)_{y})(x_{i,j}, t),\\
  &{\mathcal{L}}_1(U)_{3,i,j}=
    \rho_{i,j} G
    \approx\rho(x_{i,j},t) G
    ,\\
  &{\mathcal{L}}_2(U)_{2,i,j}\approx 
    \varepsilon(
    \frac{1}{2}(c_{y}^2)_{x}-\frac{1}{2} (c_{x}^2)_x- (c_{x}c_y)_{y})(x_{i,j}, t),\\
  &{\mathcal{L}}_2(U)_{3,i,j}\approx 
    \varepsilon(      \frac{1}{2}(c_{x}^2)_{y}-\frac{1}{2} (c_{y}^2)_y    -(c_xc_y)_x)(x_{i,j}, t),\\
  &{\mathcal{L}}_{3}(U)_{4,i,j}
    \approx\Delta(\psi'(c) - \frac{\varepsilon}{\rho}  \Delta c)  (x_{i,j}, t)
    ,\\
  &{\mathcal{L}}_4(U)_{2,i,j}\approx
    (\nu(  (v_1)_{xx}+(v_1)_{yy})+(\nu+\lambda)(  (v_1)_{xx}+(v_2)_{xy}))(x_{i,j}, t),\\
  &{\mathcal{L}}_4(U)_{3,i,j}\approx
    (\nu((v_2)_{xx}+(v_2)_{yy})  +(\nu+\lambda)(  (v_1)_{xy}+(v_2)_{yy}))(x_{i,j}, t).
\end{align*}
Here, as in the rest of this section, we drop the dependence of $U$ on $t$ to obtain
the cited spatial semidiscretization. 

The convective term ${\mathcal{C}}$ is obtained through finite differences of
numerical fluxes obtained by WENO5 reconstructions
\cite{BBMZ19a,BBMZ19b} on  Global
Lax-Friedrichs flux splittings \cite{Shu09}, which is fifth-order
accurate for finite difference schemes, based on point values.  Other schemes
for systems of hyperbolic conservation laws could be used as well, see
\cite{Toro09} and references therein.

To approximate the terms that involve derivatives of $c$ in the
conservation of momenta, we define finite difference operators for
functions on $M\times M$ grids, for
fixed $h>0$, to approximate first-order derivatives:
\begin{align*}
  &D_{x}^{1*} f_{i,j}=
    \begin{cases}
      \frac{f_{i,j}}{h}& i=1,\\
      \frac{f_{i,j}-f_{i-1,j}}{h}& 1<i<M,\\
      \frac{-f_{i-1,j}}{h}& i=M.\\
    \end{cases}
  \\
  &D_{x}^1 f_{i,j}=
    \begin{cases}
      \frac{f_{i+1,j}-f_{i,j}}{h}& i<M,\\
      0& i=M.
    \end{cases}
  \\
  &D_{x} f_{i,j}=
    \begin{cases}
      \frac{f_{i+1,j}-f_{i-1,j}}{2h}& 1<i<M,\\
      \frac{f_{i+1,j}-f_{i,j}}{h}& i=1,\\
      \frac{f_{i,j}-f_{i-1,j}}{h}& i=M.\\
    \end{cases}  
  \\
  &D^*_{x} f_{i,j}=
    \begin{cases}
      \frac{f_{i+1,j}-f_{i-1,j}}{2h}& 1<i<M,\\
      \frac{f_{i+1,j}-f_{i,j}}{2h}& i=1,\\
      \frac{f_{i,j}-f_{i-1,j}}{2h}& i=M.\\
    \end{cases}\\
  &S_{x}f_{i,j}=\begin{cases}
    f_{i+1,j} & i < M,\\
    0& i=M.
  \end{cases}    
\end{align*}
and likewise in the $y$ direction.
\begin{enumerate}
\item $D_{x}^{1*}f_{i,j}$ is a
  second-order accurate
  approximation for $f_{x}(x_{i-\frac{1}{2},j})$ when $f_{i,j}=f(x_{i,j})$ and
  $f\in {\mathcal{C}}^3$ with $f(x_{0,j})=f(x_{M,j})=0$, which is used to
  approximate pure double derivatives.
\item 
  $D_{x}^1f_{i,j}$ is a second-order accurate
  approximation for $f_{x}(x_{i+\frac{1}{2},j})$ when $f_{i,j}=f(x_{i,j})$ and
  $f\in {\mathcal{C}}^3$ with $f_{x}(x_{M+\frac{1}{2},j})=0$, which is used to
  approximate pure double derivatives. These two operators are related
  as $D_{x}^{1*}=-(D_{x}^1)^T$.
\item $D_{x}f_{i,j}$ is a second-order  accurate
  approximation for $f_{x}(x_{i,j})$ for $1<i,j<M$ and first-order
  accurate otherwise, when $f_{i,j}=f(x_{i,j})$ and
  $f\in {\mathcal{C}}^3$.
\item $D_{x}^*f_{i,j}$ is a second-order accurate
  approximation for $f_{x}(x_{i,j})$ for $1<i,j<M$ or $j=1,M$ and first-order accurate otherwise, when $f_{i,j}=f(x_{i,j})$ and
  $f\in {\mathcal{C}}^3$ with $f(x_{i,\frac{1}{2}})=f(x_{i,M+\frac{1}{2}})=0$.
\end{enumerate}
For the sake of
notation, for functions $f, g$ on $M\times M$ grids we denote
$(f*g)_{i,j}=f_{i,j}g_{i,j}$.

We consider the following second-order accurate approximations at
interior points $1<i,j<M$ and first-order accurate at the rest of the  points, for
$c_{i,j}=\frac{(\rho   c)_{i,j}}{\rho_{i,j}}\approx c(x_{i,j})$, in which the boundary conditions \eqref{eq:bc} on $c$ are taken into account:
\begin{align*}
  (c_x^2)_x(x_{i,j}) &\approx (D_{x}^{1*}(D_{x}^{1}C*D_{x}^{1}C))_{i,j},\\
  (c_y^2)_x(x_{i,j}) &\approx (D_x((D_{y}^*C*D_{y}^*C)))_{i,j},\\
  (c_xc_y)_x(x_{i,j}) &\approx \frac{1}{2} (D_x^{1*}(D_{x}^1C*(S_{x}D_y^*C+D_{y}^*C)))_{i,j},\\
  (c_y^2)_y(y_{i,j}) &\approx (D_{y}^{1*}(D_{y}^{1}C*D_{y}^{1}C))_{i,j},\\
  (c_x^2)_y(y_{i,j}) &\approx (D_y((D_{x}^*C*D_{x}^*C)))_{i,j},\\
  (c_yc_x)_y(y_{i,j}) &\approx \frac{1}{2} (D_y^{1*}(D_{y}^1C*(S_{y}D_x^*C+D_{x}^*C)))_{i,j}
\end{align*}

To approximate the terms that involve derivatives of $v$ in the
conservation of momenta, we consider the following finite difference approximation to
$(v_{k})_{xx}$, for $w=V_{k}, k=1,2,$,
which is
second-order accurate for $1<i,j<M$ and first-order accurate otherwise under the no-slip boundary
conditions on $v_{k}, k=1,2$
\begin{align*}
  &E_{x}w=
    \begin{cases}
      \frac{1}{h^2}\left(
        \frac{4}{3}w_{i+1,j}-4 w_{i,j}
      \right) & i=1,\\
      \frac{1}{h^2}\left(
        w_{i+1,j}-2 w_{i,j}+w_{i-1,j}
      \right) & 1<i< M,\\
      \frac{1}{h^2}\left(
        -4 w_{i,j}+\frac{4}{3}w_{i-1,j}
      \right) & i= M,\\
    \end{cases}
\end{align*}
The $E_{y}$ operator is defined analogously.

The approximations are:

\begin{align*}
  &(v_{k})_{xx}(x_{i,j})\approx 
    (E_{x}V_k)_{i,j},\\
  &(v_{k})_{yy}(x_{i,j})\approx 
    (E_{y}V_k)_{i,j},\\
  &(v_{k})_{xy}(x_{i,j})\approx (D_x(D_y V_k))_{i,j},
\end{align*}
which lead to the $k=2,3$ $M^2$ blocks ${\mathcal{L}}^{k}_4(U)$ of  ${\mathcal{L}}_4(U)$
\begin{align*}
  &
    \begin{bmatrix}
      {\mathcal{L}}^2_{4}(U),\\
      {\mathcal{L}}^3_{4}(U)
    \end{bmatrix}    
  =
  \begin{bmatrix}
    (2\nu+\lambda)I_{M}\otimes E+\nu E\otimes I_{M} &(\nu+\lambda)D\otimes D \\
    (\nu+\lambda)D\otimes D & \nu I_{M}\otimes E+(2\nu+\lambda) E\otimes I_{M}
  \end{bmatrix}
                              \begin{bmatrix}
                                V_1\\
                                V_2
                              \end{bmatrix},
\end{align*}
where $I_M$ is the $M\times M$ identity matrix, $\otimes$ is the Kronecker product and 
\begin{align}
  \label{eq:11b}
  &  E=
    \frac{1}{h^2}
    \begin{bmatrix}
      -4&\frac{4}{3}&0&\dots&0\\
      1 &-2&1&\dots&0\\
      \hdotsfor{5}\\
      0&\dots&1 &-2&1\\
      0&\dots&0 &\frac{4}{3}&-4
    \end{bmatrix},
                              \quad      
                              D=
                              \frac{1}{h}
                              \begin{bmatrix}
                                -1&1&0&\dots&0\\
                                -\frac{1}{2} &0&\frac{1}{2}&\dots&0\\
                                \hdotsfor{5}\\
                                0&\dots&-\frac{1}{2} &0&\frac{1}{2}\\
                                0&\dots&0 &-1&1
                              \end{bmatrix}.
\end{align}
The matrices in \eqref{eq:11b} fail to be symmetric due to the boundary
conditions. This could be circumvented with staggered grids for the
velocity, but we do not consider this possibility in this
paper. 

The term $(\psi'(c) )_{xx}=(\psi''(c)c_x)_x$ needs special care since
it is not negatively definite due to $\psi''(c)=3c^2-1$ changing sign
in $(-1,1)$. Following \cite{Vollmayr-Lee-Rutenberg2003} we 
consider the splitting $\psi'=\phi_{+}+\phi_{-}$\begin{align*}
  &\phi_{-}=c^3-3c, \phi_{+}=2c, \\
  &\phi_{-}'(c)=3(c^2-1)\leq 0,
    \phi_{+}'(c)=2> 0\,\forall c\in[-1,1].
\end{align*}        

For $\chi=\phi_{\pm}$, 
taking into account that $(\chi(c))_{z}=\chi'(c)c_z$, $z=x,y$, so $\chi(c)$
satisfies  Neumann boundary conditions, 
we have the following second-order accurate approximations:
\begin{equation}\label{eq:222}
  \begin{aligned}
    &(\chi(c))_{xx}(x_{i,j})
    \\
    &\approx 
    \begin{cases}
      \frac{(\chi'(c_{i+1,j})+\chi'(c_{i,j}))(c_{i+1,j}-c_{i,j})}{2h^2}&i=1,\\
      \frac{(\chi'(c_{i+1,j})+\chi'(c_{i,j}))(c_{i+1,j}-c_{i,j})-(\chi'(c_{i,j})+\chi'(c_{i-1,j}))(c_{i,j}-c_{i-1,j})}{2h^2}&1<i<M,\\
      \frac{-(\chi'(c_{i,j})+\chi'(c_{i-1,j}))(c_{i,j}-c_{i-1,j})}{2h^2}&i=M.\\
    \end{cases}
    \\
    &(\chi(c))_{yy}(x_{i,j})\\
    &\approx 
    \begin{cases}
      \frac{(\chi'(c_{i,j+1})+\chi'(c_{i,j}))(c_{i,j+1}-c_{i,j})}{2h^2}&j=1,\\
      \frac{(\chi'(c_{i,j+1})+\chi'(c_{i,j}))(c_{i,j+1}-c_{i,j})-(\chi'(c_{i,j})+\chi'(c_{i,j-1}))(c_{i,j}-c_{i,j-1})}{2h^2}&1<j<M,\\
      \frac{-(\chi'(c_{i,j})+\chi'(c_{i,j-1}))(c_{i,j}-c_{i,j-1})}{2h^2}&j=M.\\
    \end{cases}
  \end{aligned}
\end{equation}
These yield approximations
\begin{align}\label{eq:223}
  (\mathcal{M}_{\pm}(C)C)_{i,j} \approx
  \Delta(\phi_{\pm}(c))(x_{i,j}, t),
\end{align}  
where we denote by
$\mathcal{M}_{\pm}(C)$  the tensor built from the values
of $\chi_{\pm}$ that appear in \eqref{eq:222}.

It can be seen that the boundary conditions $\nabla c(x, y, t)\cdot
\boldsymbol{n}(x, y)=\nabla \mu(x,y,t) \cdot \boldsymbol{n}(x,y)=0$ are equivalent to 
$\nabla c(x, y, t)\cdot
\boldsymbol{n}(x,y)=\nabla \xi(x,y,t) \cdot \boldsymbol{n}(x,y)=0$
for $\xi=\frac{1}{\rho}\Delta c$.

If $f(x,y)$ satisfies the Neumann boundary condition $\nabla
f(x,y)\cdot \boldsymbol{n}(x,y)=0$ and $f_{i,j}=f(x_{i,j})$ then
$(\Delta_h f)_{i,j}$, $\Delta_h=\Delta_{x,h}+\Delta_{y,h}$,  is a second-order accurate
approximation of $\Delta f(x_{i,j})$ for $f\in{\mathcal{C}}^4$, where
\begin{align*}
  &
    \Delta_{x,h} f_{i,j}=
    \begin{cases}
      \frac{f_{i+1,j}-f_{i,j}}{h^2} & i=1,\\
      \frac{f_{i+1,j}-2f_{i,j}+f_{i-1,j}}{h^2} & 1<i<M,\\
      \frac{f_{i-1,j}-f_{i,j}}{h^2} & i=M.
    \end{cases}\\
  &
    \Delta_{y,h} f_{i,j}=
    \begin{cases}
      \frac{f_{i,j+1}-f_{i,j}}{h^2} & j=1,\\
      \frac{f_{i,j+1}-2f_{i,j}+f_{i,j-1}}{h^2} & 1<j<M,\\
      \frac{f_{i,j-1}-f_{i,j}}{h^2} & j=M.
    \end{cases}
\end{align*}
Therefore we get 
\begin{align*}
  &\Delta\xi(x_{i,j})\approx \Big(\Delta_h\big(D(\varrho)^{-1}\Delta_h
    c\big)\Big)_{i,j},
\end{align*}    
which, together with \eqref{eq:223}, yields the approximation
\begin{align}\label{eq:224}
  &\Delta \mu(x_{i,j})\approx
    {\mathcal{L}}^4_3(U)_{i,j}=\Big(\mathcal{M}_{+}(C)C+\mathcal{M}_{-}(C)C-\varepsilon\Delta_h\big(
    D(\varrho)^{-1}\Delta_h C\big)\Big)_{i,j},
\end{align}
where $D(v), v\in{\mathbb{R}}^{M\times M},$ is the diagonal operator on $M\times
M$ matrices given by
\begin{align*}
  (D(v) w)_{i,j}=v_{i,j}w_{i,j}, i,j =1,\dots,M.
\end{align*}

The numerical schemes for the 1D case are obtained in a
straightforward manner.

\subsection{IMEX schemes}
A Runge-Kutta solver, as the Explicit Euler method
\begin{align}\label{eq:ee}
  U^{n+1}&=U^{n}+\Delta t \mathcal{L}(U^{n}),
\end{align}
applied to obtain a fully discrete scheme
would  require $\Delta t \propto \Delta x^4$ for stability, which would yield a prohibitively expensive numerical scheme. 

Instead, we use the technique of doubling variables and partitioned Runge-Kutta
schemes \cite{BBMRV15,PR05} to obtain Linearly IMplicit EXplicit schemes.
We denote  the variables that are to be
treated explicitly with a tilde and define
\begin{equation*}
  \widetilde{\mathcal{L}}(\widetilde U, U)=
  {\mathcal{C}}(\widetilde U)+{\mathcal{L}}_1(U)+{\mathcal{L}}_2(U)+
  \widetilde{{\mathcal{L}}}_3(\widetilde U, U)+{\mathcal{L}}_4(U),
\end{equation*}
where 
\begin{align*}
  \widetilde U=
  \begin{bmatrix}
    \widetilde \varrho\\
    \widetilde \varrho * \widetilde V_1\\
    \widetilde \varrho * \widetilde V_2\\
    \widetilde \varrho * \widetilde C
  \end{bmatrix}    ,
  \quad
  U=
  \begin{bmatrix}
    \varrho\\
    \varrho * V_1\\
    \varrho * V_2\\
    \varrho * C
  \end{bmatrix},\quad
  \widetilde{\mathcal{L}}_{3}(\widetilde U, U)
  =\begin{bmatrix}
    0\\
    0\\
    0\\
    \widetilde{\mathcal{L}}_{3}^4(\widetilde U, U)
  \end{bmatrix}    
\end{align*}
with 
\begin{equation}\label{eq:3}
  \widetilde{\mathcal{L}}_{3}^4(\widetilde U, U)=
  \mathcal{M}_{+}(\widetilde{C})C+\mathcal{M}_{-}(\widetilde{C})\widetilde{C}-\varepsilon\Delta_h\big(
  D(\varrho)^{-1}\Delta_h C\big).
\end{equation}

From \eqref{eq:224} and \eqref{eq:2}, the requirement  $\widetilde{\mathcal{L}}(U, U)={\mathcal{L}}(U)$ is met.

Now we have the IVP
\begin{equation}\label{eq:225}
  \begin{aligned}
    \widetilde U'&=\widetilde {\mathcal{L}}(\widetilde U, U)\\
    U'&=\widetilde {\mathcal{L}}(\widetilde U, U)\\
    \widetilde U(0)&=U(0)=U_0\\
  \end{aligned}
\end{equation}
is equivalent to 
\begin{align*}
  U'&= {\mathcal{L}}(U)\\
  U(0)&=U_0.
\end{align*}

A partitioned Runge-Kutta scheme, in
which there are two different $s$ stages Butcher tableaus, one explicit
and one (diagonally) implicit
\begin{equation*}
  \begin{array}{c|c}
    \widetilde \gamma & \widetilde \alpha\\
    \hline
                      &\widetilde \beta^T
  \end{array},\quad \widetilde\alpha_{i,j}=0, j\geq i
  \qquad
  \begin{array}{c|c}
    \gamma &  \alpha\\
    \hline
           &\beta^T
  \end{array},\quad \alpha_{i,j}=0, j> i,
\end{equation*}
can be applied to \eqref{eq:225}.
It can be seen that if both Butcher tableaus yield second-order accurate schemes
and $\beta=\widetilde{\beta}$, then the resulting partitioned
Runge-Kutta scheme is second-order accurate.

This results in the recursion, for $i=1,\dots,s$:
\begin{align*}
  \widetilde U^{(i)}&=\widetilde U^{n}+\Delta t \sum_{j<i}\widetilde \alpha_{i,j}
                      \widetilde {\mathcal{L}}(\widetilde U^{(j)}, U^{(j)}),\\
  U^{(i)}&=U^{n}+\Delta t \sum_{j< i} \alpha_{i,j}
           \widetilde {\mathcal{L}}(\widetilde U^{(j)}, U^{(j)})
           +\Delta t \alpha_{i,i}
           \widetilde {\mathcal{L}}(\widetilde U^{(i)}, U^{(i)}),\\
  \widetilde U^{n+1}&=\widetilde U^{n}+\Delta t\sum_{j=1}^{s} \widetilde \beta_{j}\widetilde{\mathcal{L}}(\widetilde U^{(j)}, U^{(j)}),\\
  U^{n+1}&=U^{n}+\Delta t\sum_{j=1}^{s} \beta_{j}\widetilde{\mathcal{L}}(\widetilde U^{(j)}, U^{(j)})=\widetilde U^{n+1},
\end{align*}
since $\widetilde \beta_{j}=\beta_j$, $\forall j$ and $U^{n}=\widetilde
U^{n}$,  and there is no need
of doubling variables, which we henceforth assume.

The definitive recursion is the following:
\begin{align*}
  \widetilde U^{(i)}&= U^{n}+\Delta t \sum_{j<i}\widetilde \alpha_{i,j}  {\mathcal{K}}_j,\\
  U^{(i)}&=U^{n}+\Delta t \sum_{j< i} \alpha_{i,j}
           {\mathcal{K}}_j  +\Delta t \alpha_{i,i}
           \widetilde {\mathcal{L}}(\widetilde U^{(i)}, U^{(i)}),\\
  U^{n+1}&=U^{n}+\Delta t\sum_{j=1}^{s} \beta_{j} {\mathcal{K}}_j,
\end{align*}
where
\begin{align*}
  {\mathcal{K}}_j=\widetilde {\mathcal{L}}(\widetilde U^{(j)}, U^{(j)}).
\end{align*}

We consider in this paper Stiffly Accurate Runge-Kutta solvers, i.e.,
the last row of the $\alpha$ matrix coincides with
$\beta^T$. Specifically, we consider the following Butcher tableaus
\begin{align*}
  \text{EE-IE }&
                    \begin{array}{c|c}
                      0& 0\\
                      \hline
                       & 1
                    \end{array}    
                       && 
                          \begin{array}{c|c}
                            1& 1\\
                            \hline
                             & 1
                          \end{array}    
  \\
  \text{*-DIRKSA}     &\begin{array}{c|cc}
                            0&0& 0\\
                            1+s&    1+s& 0\\
                            \hline
                             &    s&     1-s
                          \end{array}    
                       && 
                          \begin{array}{c|cc}
                            1-s&1-s& 0\\
                            1&    s&     1-s\\
                            \hline
                               &    s&     1-s
                          \end{array},
                                       \quad s=\frac{1}{\sqrt{2}}.
\end{align*}
The DIRKSA scheme is the 
only 2-stages second order stiffly accurate DIRK method with $\alpha_{ij}\geq 0$.  

The order for EE-IE is 1 and it is 2 for *-DIRKSA.

\subsection{Systems solutions}

One needs to solve
\begin{align}\label{eq:4}
  U^{(i)}
  =U^{n}+\Delta t \sum_{j< i} \alpha_{i,j}{\mathcal{K}}_j
  +\Delta t \alpha_{i,i}
  \widetilde {\mathcal{L}}(\widetilde U^{(i)}, U^{(i)}),
\end{align}
for $U^{(i)}$, where
\begin{align*}
  U^{n}=
  \begin{bmatrix}
    \varrho^{n}\\
    M_1^{n}\\
    M_2^{n}\\
    Q^{n}
  \end{bmatrix},
  \quad
  U^{(i)}=
  \begin{bmatrix}
    \varrho^{(i)}\\
    M_1^{(i)}\\
    M_2^{(i)}\\
    Q^{(i)}
  \end{bmatrix}=
  \begin{bmatrix}
    \varrho^{(i)}\\
    \varrho^{(i)}*V_1^{(i)}\\
    \varrho^{(i)}*V_2^{(i)}\\
    \varrho^{(i)}* C^{(i)}
  \end{bmatrix}.    
\end{align*}  
As we shall see, although ${\mathcal{L}}_2, \widetilde{{\mathcal{L}}}_3, {\mathcal{L}}_4$ are not linear,
only  linear systems for $V_k^{(i)},
C^{(i)}$  have to be solved. 

For the first variable, with block superscript notation for the
operators and ${\mathcal{K}}$ variables, we get:
\begin{align*}
  &\varrho^{(i)}= 
    \varrho^{n}+\Delta t \sum_{j< i} \alpha_{i,j}
    {\mathcal{K}}^1_j+\Delta t \alpha_{i,i}
    \widetilde {\mathcal{C}}^1(\widetilde U^{(i)}),
\end{align*}
so $\varrho^{(i)}$ is explicitly computable.

For the fourth variable $Q^{(i)}$, since $\varrho^{(i)}$ is already known, this
system can be cast for the 
$C^{(i)}$  variables:
\begin{align*}
  Q^{(i)}
  &=Q^{n}+\Delta t \sum_{j< i} \alpha_{i,j}
    {\mathcal{K}}^4_j \\
  &+ \Delta t \alpha_{i,i}
    \Big(
    \widetilde {\mathcal{C}}^4(\widetilde U^{(i)})
    +
    \mathcal{M}_{+}(\widetilde{C}^{(i)}) C^{(i)}+
    \mathcal{M}_{-}(\widetilde{C}^{(i)}) \widetilde{C}^{(i)}-\varepsilon\Delta_h
    D(\varrho^{(i)})^{-1}\Delta_hC^{(i)}
    \Big),
\end{align*}
which is equivalent to
\begin{equation}\label{eq:226}
  \begin{aligned}
    &\left(D(\varrho^{(i)})-\Delta t \alpha_{i,i}
      \mathcal{M}_{+}(\widetilde{C}^{(i)}) +\Delta t \alpha_{i,i}\varepsilon\Delta_h
      D(\varrho^{(i)})^{-1}\Delta_h\right)C^{(i)}\\
    &=Q^{n}+\Delta t \sum_{j< i} \alpha_{i,j}
    {\mathcal{K}}^4_j
    + \Delta t \alpha_{i,i}
    \Big(
    \widetilde {\mathcal{C}}^4(\widetilde U^{(i)})
    +\mathcal{M}_{-}(\widetilde{C}^{(i)}) \widetilde{C}^{(i)}
    \Big).
  \end{aligned}
\end{equation}
If $\varrho^{i}_k>0\,\forall k$, then the matrix of this system is
symmetric and positive definite, for is  the sum of a diagonal positive matrix and two symmetric and positive
semidefinite matrices.

For the second and third variables, one needs to solve:
\begin{align*}
  &\begin{bmatrix}
    \varrho^{(i)}*V_1^{(i)}\\
    \varrho^{(i)}*V_2^{(i)}
  \end{bmatrix}    
  =\begin{bmatrix}
    M_1^{n}\\
    M_2^{n}
  \end{bmatrix}
  +\Delta t \sum_{j< i} \alpha_{i,j}
  \begin{bmatrix}
    {\mathcal{K}}^2_j\\
    {\mathcal{K}}^3_j
  \end{bmatrix}    
  + \Delta t \alpha_{i,i}
  \begin{bmatrix}
    {\mathcal{C}}^2(\widetilde U^{(i)})+
    {\mathcal{L}}_1^2(U^{(i)})+
    {\mathcal{L}}_2^2(U^{(i)})
    \\
    {\mathcal{C}}^3(\widetilde U^{(i)})+
    {\mathcal{L}}_1^3(U^{(i)})+
    {\mathcal{L}}_2^3(U^{(i)})
  \end{bmatrix}
  \\
  &+\Delta t \alpha_{i,i}
    \begin{bmatrix}
      (2\nu+\lambda)I_{M}\otimes E+\nu E\otimes I_{M} &(\nu+\lambda)D\otimes D \\
      (\nu+\lambda)D\otimes D & \nu I_{M}\otimes E+(2\nu+\lambda) E\otimes I_{M}
    \end{bmatrix}
                                \begin{bmatrix}
                                  v_{1}^{(i)}\\
                                  v_{2}^{(i)}
                                \end{bmatrix}.
\end{align*}
Since      ${\mathcal{L}}_1^j(U^{(i)})$ and      ${\mathcal{L}}_2^j(U^{(i)}), j=2,3,$ do not depend
on $V_1^{(i)},V_2^{(i)}$, they can be computed from  previous steps and   there only
remains the following equation to be solved:
\begin{equation}\label{eq:227}
  \begin{aligned}
    &\left(\begin{bmatrix}
        D(\varrho^{(i)})&0\\
        0&D(\varrho^{(i)})
      \end{bmatrix}
    \right.
    \\
    &\left.  -\Delta t \alpha_{i,i}
      \begin{bmatrix}
        (2\nu+\lambda)I_{M}\otimes E+\nu E\otimes I_{M} &(\nu+\lambda)D\otimes D \\
        (\nu+\lambda)D\otimes D & \nu I_{M}\otimes E+(2\nu+\lambda) E\otimes I_{M}
      \end{bmatrix}
    \right)
    \begin{bmatrix}
      V_{1}^{(i)}\\
      V_{2}^{(i)}
    \end{bmatrix}
    \\
    &=\begin{bmatrix}
      M_1^{n}\\
      M_2^{n}
    \end{bmatrix}
    +\Delta t \sum_{j< i} \alpha_{i,j}
    \begin{bmatrix}
      {\mathcal{K}}^2_j\\
      {\mathcal{K}}^3_j
    \end{bmatrix}    
    + \Delta t \alpha_{i,i}
    \begin{bmatrix}
      {\mathcal{C}}^2(\widetilde U^{(i)})+
      {\mathcal{L}}_1^2(U^{(i)})+
      {\mathcal{L}}_2^2(U^{(i)})
      \\
      {\mathcal{C}}^3(\widetilde U^{(i)})+
      {\mathcal{L}}_1^3(U^{(i)})+
      {\mathcal{L}}_2^3(U^{(i)})
    \end{bmatrix}.
  \end{aligned}
\end{equation}
If $\varrho^{(i)}_k>0\,\forall k$, then  the matrix of this system
should be close to  symmetric and positive definite, since the matrix
\begin{align*}
  -\begin{bmatrix}
    (2\nu+\lambda)I_{M}\otimes E+\nu E\otimes I_{M} &(\nu+\lambda)D\otimes D \\
    (\nu+\lambda)D\otimes D & \nu I_{M}\otimes E+(2\nu+\lambda) E\otimes I_{M}
  \end{bmatrix}
\end{align*}                              
is a discretization  of the  self-adjoint elliptic operator
\begin{align*}
  -\Big((\nu+\lambda)\nabla \mathop{\text{div}} \boldsymbol{v}+\nu \Delta    \boldsymbol{v}\Big),
\end{align*}
under the boundary conditions \eqref{eq:bc}.

\subsection{Linear solvers}\label{ss:solver}
We have used the multigrid V-cycle algorithm with 4 pre- and post-
Gauss-Seidel smoothings and direct solution when the size of the
projected systems is $\leq 4$ (see \cite{BrandtLivne2011}) for the solution
of systems \eqref{eq:226}  and \eqref{eq:227}. 

The matrix in system \eqref{eq:226} is symmetric and positive, so
there is the possibility of using the conjugate gradient method, using 
 the following approximation 
\begin{align*}
  A&=D(\varrho^{(i)})-\Delta t \alpha_{i,i}
    \mathcal{M}_{+}(\widetilde{C}^{(i)}) +\Delta t \alpha_{i,i}\varepsilon\Delta_h
    D(\varrho^{(i)})^{-1}\Delta_h
  \\
  &\approx B=\mu_1 I_{M^2}+\Delta t \alpha_{i,i} \mu_2\Delta_h+\Delta t
    \alpha_{i,i}\varepsilon\mu_3\Delta_h^2,\\
  &\mu_1=\text{mean}(\varrho^{(i)}),\quad
    \mu_2=\text{mean}(\phi'_+(\widetilde{C}^{(i)}))=2,\quad \mu_3=\text{mean}(1/\varrho^{(i)}),
\end{align*}
as preconditioner, for $\Delta_h$ can be efficiently diagonalized by
discrete cosine transforms.

We next analyze the condition number of the preconditioned matrix:
\begin{align*}
  \frac{\max\limits_{z\neq 0}\frac{z^TAz}{z^TBz}}{\min\limits_{z\neq 0}\frac{z^TAz}{z^TBz}}.
\end{align*}

We drop the superindex $^{(i)}$ for simplicity. The matrix $\mathcal{M}_{+}(\widetilde{C})$ can be expressed as
\begin{align*}
  &\mathcal{M}_{+}(\widetilde{C}) = -(I_M\otimes
  D_1^T)D(\lambda^x)(I_M\otimes D_1)
    -(D_1^T\otimes I_M)D(\lambda^y)(D_1\otimes I_M)\\
  &D_1=
    \frac{1}{h}
    \begin{bmatrix}
      -1& 1&0&\dots&0\\
      0&-1& 1&\dots&0\\
      \hdotsfor{5}\\
      0&\dots&0&-1& 1\\
      0&\dots&0&0& 0
    \end{bmatrix}\in{\mathbb{R}}^{M\times M},
  \\
  &
    \lambda^x_{i+M(j-1)}=
    \begin{cases}
      \frac{1}{2}(\phi_+'(C_{i+M(j-1)})+\phi_+'(C_{i+1+M(j-1)}))
      & i<M\\
      0&i=M
    \end{cases}                     
  \\
  &
         \lambda^y_{i+M(j-1)}=
         \begin{cases}
           \frac{1}{2}(\phi_+'(C_{i+M(j-1)})+\phi_+'(C_{i+M(j)}))
           & j<M\\
           0&j=M;
         \end{cases}.
\end{align*}
Since $\phi_+'=2$ and 
\begin{align*}
  -(I_M\otimes  D_1^T)(I_M\otimes D_1)
    -(D_1^T\otimes I_M)(D_1\otimes I_M)=\Delta_h,
\end{align*}
we get the following:
\begin{align*}
  z^TAz&=\sum_{k=1}^{M^2} \varrho_{k}z_k^2+2\Delta t \alpha_{i,i}
    (\sum_{k=1}^{M^2} ((I_M\otimes D_1)z)_k^2+((D_1\otimes I_M)z)_k^2
    ) \\
  &+\Delta t \alpha_{i,i}\varepsilon
\Delta_h \sum_{k=1}^{M^2} \frac{1}{\varrho_k} (\Delta_hz)_k^2\\
       z^TBz&=\mu_1\sum_{k=1}^{M^2} z_k^2+2\Delta t \alpha_{i,i}
    (\sum_{k=1}^{M^2} ((I_M\otimes D_1)z)_k^2+((D_1\otimes I_M)z)_k^2
    ) \\
  &+\mu_3\Delta t \alpha_{i,i}\varepsilon
\Delta_h \sum_{k=1}^{M^2} (\Delta_hz)_k^2.
\end{align*}
Therefore, for $0\neq z\in{\mathbb{R}}^{M^2}$:
\begin{align*}
  \min\left(\frac{\min \varrho_j}{\mu_1}, \frac{\min
  \frac{1}{\varrho_j}}{\mu_3}\right)
\leq \frac{z^TAz}{z^TBz}  \leq 
  \max\left(\frac{\max \varrho_j}{\mu_1}, \frac{\max
  \frac{1}{\varrho_j}}{\mu_3}\right)
\end{align*}
therefore the condition of the preconditioned matrix is bounded above
by
\begin{align*}
\frac{\max\left(\frac{\max \varrho_j}{\mu_1}, \frac{\max
  \frac{1}{\varrho_j}}{\mu_3}\right)}{  \min\left(\frac{\min \varrho_j}{\mu_1}, \frac{\min
  \frac{1}{\varrho_j}}{\mu_3}\right)},
\end{align*}
which is close to 1 if $\varrho$ is nearly constant. Therefore, it
is expected to be a good preconditioner in this case.

\subsection{Time-step selection}\label{ss:tss}
The time-step  stability restrictions of  the purely convective
  subsystem is 
  \begin{align}\label{eq:cfl}
    \Delta t =\text{CFL} \cdot \text{cs}\cdot \Delta x,
  \end{align}  
  where CFL is a constant and the maximum of the characteristic speeds,
  $\text{cs}$, is computed, at each Runge-Kutta step, as
  \begin{align*}
    \text{cs}=\max_{k=1,2,j=1,\dots,M^2}|V_{k,j}^{(i)}|+\sqrt{\gamma(\varrho_{j}^{(i)})^{\gamma-1}}.
  \end{align*}

  The scheme is not ensured to be bound preserving, i.e., it might
  happen that density might become negative or the $c$-variable be
  outside $[-1,1]$. Purely convective models might develop vacuum
  regions and coping with this possibility is certainly challenging.

  In our case, there is no guarantee that the solution of
  \eqref{eq:226} be in $[-1,1]$. We have used in our simulation the
  strategy of decreasing $\Delta t$ when $|c|$ reaches some threshold
  (1.5 in the experiments) and increasing it until the maximum otherwise.

\section{Numerical experiments}\label{s:numerics}
The objectives  of the experiments in this section are the following:
\begin{enumerate}
\item Showing that the order of the global errors in some
  experiments coincides with the expected design order of the scheme
  used to obtain them.

\item Showing that some IMEX schemes can perform time steps $\Delta
  t$ with the same stability restrictions as the purely convective
  subsystem, see \eqref{eq:cfl}.
  \item Testing the behavior of different issues for the algorithms, such as conservation, number of iterations for the linear solvers, etc.
\end{enumerate}

In all numerical experiments, the adiabatic constant $\gamma$ has been
set to $5/3$.

All the results have been obtained with a C++ implementation, using
the GNU C++ compiler with optimizations -O3 and running in a single
core of an AMD EPYC 7282 3.0 GHz CPU. The matrices of systems
\eqref{eq:226} and \eqref{eq:227} are stored by
diagonals.

\subsection{Stability test}
We consider the following initial condition for a one-dimensional test:
($c_0$ in unstable region $(-\frac{1}{\sqrt 3}, \frac{1}{\sqrt 3})$)
\begin{align*}
  &\rho_0(x)=0.1\cos(2 \pi x)+1.25\\
  &v_{0}(x) =  \sin(\pi x) \\
  &c_0(x)=0.1\cos( \pi x) 
\end{align*}
with  parameters  $G=-10$, $\nu_*=2\nu+\lambda=2$, $\varepsilon=10^{-4}$.
The Explicit Euler scheme \eqref{eq:ee} blows up for $M=8000$,  and $\Delta t=\Delta x^3$ for
$t\approx  6\cdot 10^{-11}$, thus indicating that $\Delta t$
should be proportional to $\Delta x^4$ for stable simulations.

The EE-IE and  *-DIRKSA blow up for $M=100$ and $\Delta
t$ computed by \eqref{eq:cfl} for $\text{CFL}=1.1$ for
$t\approx  10^{-1}$, whereas they do not  for  $\text{CFL}=1$ and $M=10000$.

\subsection{Order test}
This test aims to the assessment that the *-DIRKSA method achieves
second-order accuracy in the global errors. For this purpose, we add a
forcing term to the equations 
so that the solution is prescribed. Specifically, the solution in this
case is 
\begin{align*}
  \rho(x,y,t)&=\frac{\cos\left(2\,\pi \,x\right)\,\cos\left(\pi \,y\right)\,\left(t+1\right)}{10}+\frac{5}{4},\\
  v_1(x,y,t)&= -\sin\left(\pi \,x\right)\,\sin\left(\pi \,y\right)\,\left(2\,t^2-1\right),\\
  v_2(x,y,t)&= \sin\left(\pi \,x\right)\,\sin\left(2\,\pi \,y\right)\,\left(t^2+1\right),\\
  c(x,y,t)&= \frac{3}{4}-\frac{\cos\left(\pi \,x\right)\,\cos\left(\pi \,y\right)\,\left(t-1\right)}{10}.
\end{align*}  
Notice that these functions satisfy the boundary conditions
\eqref{eq:bc}.

The parameters that have been used are the following:
\begin{align*}
  \nu=1, \quad \lambda=10^{-1},\quad \varepsilon=10^{-4},\quad G=-10.
\end{align*}

For these tests, we have used $\Delta t$ given by
\eqref{eq:cfl} with CFL=0.4.

For $M\times M$ grids, with $M=2^{l}, l=3,\dots,8$, the global errors
for the approximations $u_{k,i,j}^n$ obtained by the *-DIRKSA method
for $t_n=T=0.01$,  are computed as 
\begin{align*}
  e_{M}=\frac{1}{M^2}\sum_{k=1}^{4}\sum_{i,j=1}^{M}|u_{k,i,j}^n-u_{k}(x_{i,j},
  T)|,
\end{align*}
and are displayed in Table \ref{tbl:1}, where it can be observed the
convergence of the quotients $e_{M}/e_{2M}$ towards 4. The 
analogous experiment is performed for the EE-EI scheme  resulting in quotients $e_{M}/e_{2M}$ that decrease away from 4.

\begin{table}[htb]
  \begin{center}
    \begin{tabular}{|c|c|}
      \hline
    *-DIRKSA & EE-EI\\\hline
\begin{tabular}{c|c|c}
  $M$&$e_M$ & $e_M/e_{2M}$\\\hline
  8  &  1.9828e-02  & 4.62\\
  16 &  4.2964e-03  & 4.12\\
  32 &  1.0422e-03  & 3.92\\
  64 &  2.6617e-04  & 3.93\\
  128&  6.7802e-05  & 3.95\\
  256&  1.7148e-05  & $-$
\end{tabular}
            &
\begin{tabular}{c|c|c}
  $M$&$e_M$ & $e_M/e_{2M}$\\\hline
  8  & 1.4989e-02   & 4.76\\
  16 & 3.1522e-03   & 3.38\\
  32 & 9.3289e-04   & 3.40\\
  64 & 2.7460e-04   & 3.31\\
  128& 8.2857e-05   & 3.02\\
  256& 2.7457e-05   & $-$
\end{tabular}
                      \\\hline
  \end{tabular}                      
\end{center}
\caption{Computed orders of convergence of global errors of *-DIRKSA
  and EE-EI IMEX schemes for the test with a forced solution.}
\label{tbl:1}
\end{table}

\subsection{Two-dimensional tests.}
For the following two-dimensional tests we have used $\Delta t$ given by
\eqref{eq:cfl} with CFL=0.4, which is a safe setup for simulations
with the corresponding explicit schemes for the convective part only
(isentropic Euler equations).

We consider the following tests:
\begin{itemize}
\item Test 1 ($c_0$ in unstable region $(-\frac{1}{\sqrt 3}, \frac{1}{\sqrt 3})$):
  \begin{align*}
    &\rho_0(x, y)=0.1\cos(2 \pi x) \cos(\pi y)+1.25\\
    &\boldsymbol{v}_{0}(x, y) =(  \sin(\pi x) \sin(\pi y),   \sin(\pi x)
      \sin(2\pi y))\\
    &c_0(x, y)=0.1\cos( \pi x) \cos(\pi y)
  \end{align*}
\item Test 2 ($c_0$ in stable region):
  \begin{align*}
    &\rho_0(x, y)=0.1\cos(2 \pi x) \cos(\pi y)+1.25\\
    &\boldsymbol{v}_{0}(x, y) =(  \sin(\pi x) \sin(\pi y),   \sin(\pi x)
      \sin(2\pi y))\\
    &c_0(x, y)=0.75+0.1\cos( \pi x) \cos(\pi y)
  \end{align*}
\item Test 3: $\rho=1, \boldsymbol{v}_0=0,$ $c_0$ uniform random sample of $0$ mean
  and $10^{-10}$ standard deviation.
  
\end{itemize}

Notice that these functions satisfy the boundary conditions
\eqref{eq:bc} (for test 3, almost within roundoff error).

\begin{figure}[htb]
\begin{center}
\includegraphics[height=4.5cm]{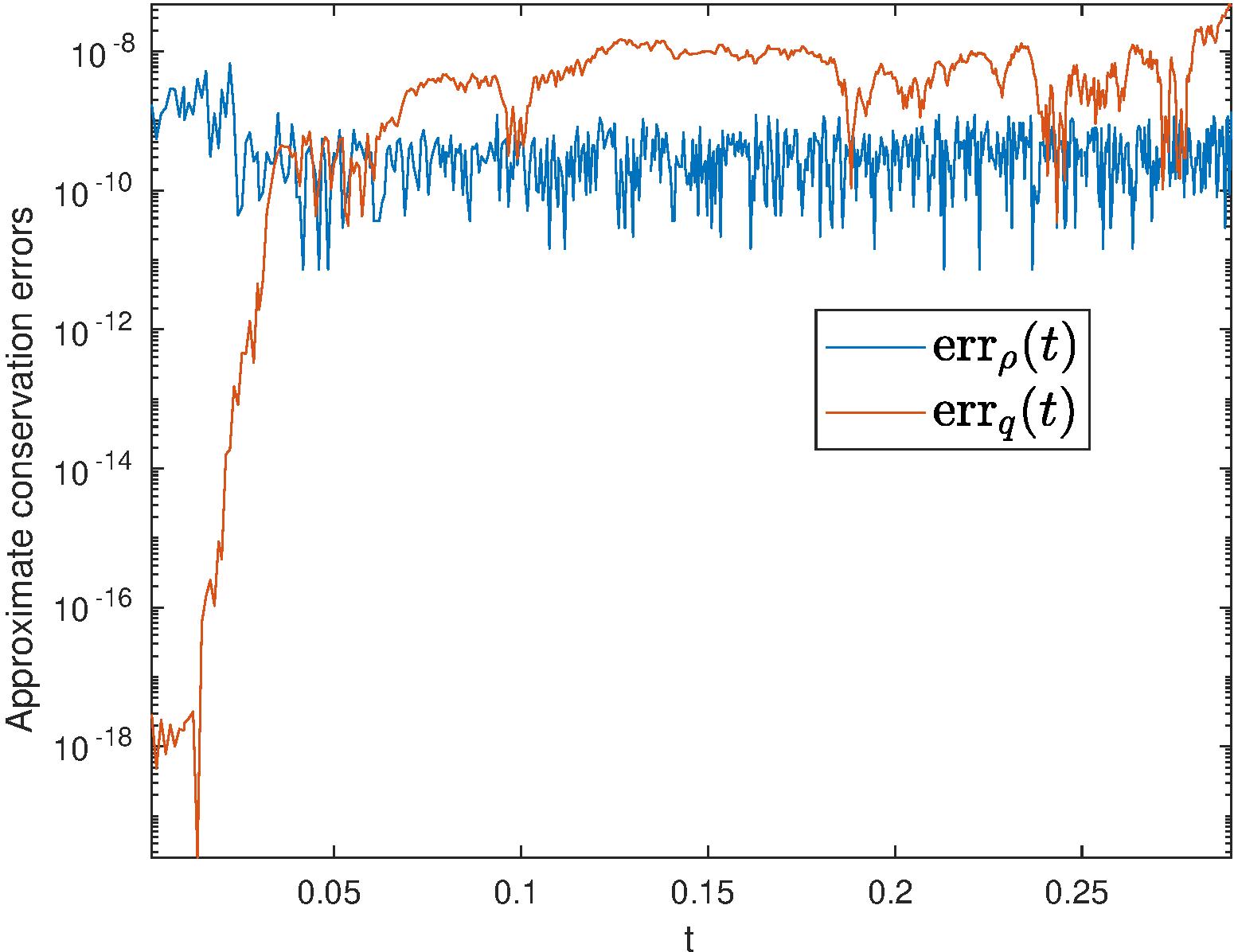}
\end{center}
\caption{Conservation errors for test 3.}
\label{fig:cons}
\end{figure}

In Figure \ref{fig:cons}, we show the time evolution of the conservation errors for $\rho$ and $q=\rho c$ for Test 3, and parameters $\nu=10^{-3}, \lambda=10^{-4}, \nu=10^{-4}, M=256$, with multigrid with relative decrease of residual of $10^{-12}$ as stopping criterion. Specifically, we approximate 
\begin{align*}
\int_{\Omega}\rho(x, t_n)dx-
  \int_{\Omega}\rho(x, 0)dx
  &\approx \text{err}_{\rho}(t_n)=\sum_{i,j=1}^{M}\varrho_{i,j}^{n}-
    \sum_{i,j=1}^{M}\varrho_{i,j}^{0}
\\
  \int_{\Omega}\rho(x, t_n)dx-
  \int_{\Omega}\rho(x, 0)dx
  &\approx \text{err}_{q}(t_n)=\sum_{i,j=1}^{M}Q_{i,j}^{n}-
    \sum_{i,j=1}^{M}Q_{i,j}^{0}
\end{align*}

\begin{figure}[htb]
\begin{center}
\includegraphics[width=0.4\hsize]{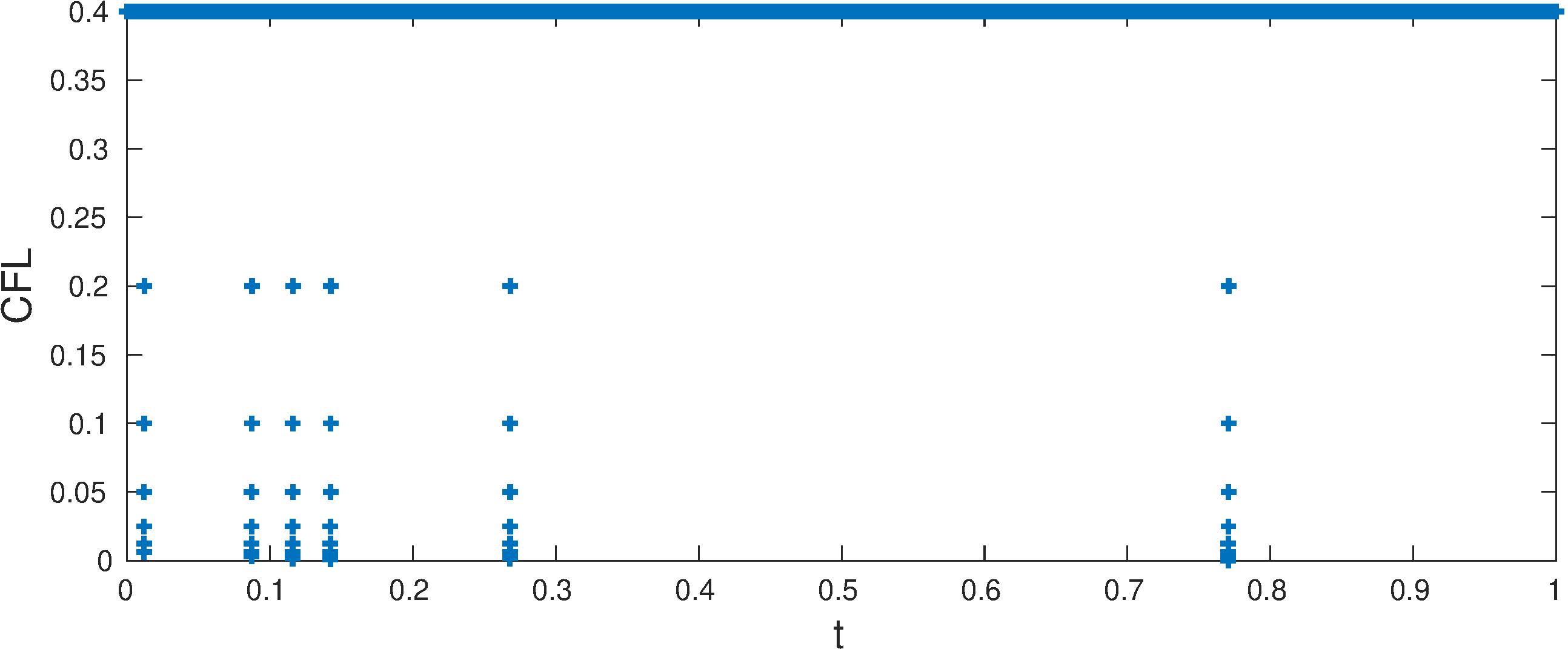}
\includegraphics[width=0.4\hsize]{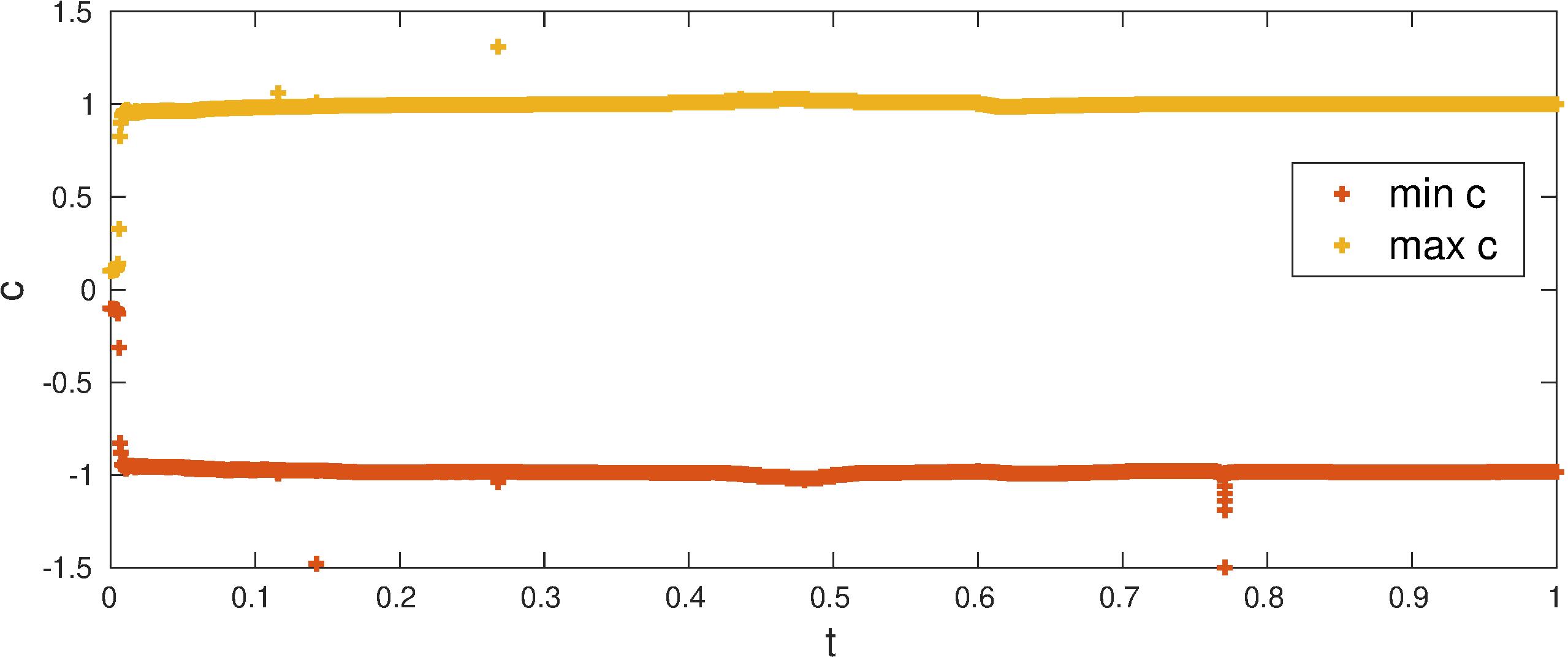}
\end{center}
\caption{Test 1. Left: Time evolution of the CFL parameter; Right:
  Time evolution of $\min c$, $\max c$.}
\label{fig:tss}
\end{figure}

In Figure \ref{fig:tss} we show the time evolution of the CFL
parameter, according to subsection \ref{ss:tss}, for a specially
challenging case for Test 1, with
parameters $\nu=10^{-3},\lambda=10^{-4}, \varepsilon=10^{-5}, M=256$.

\begin{figure}[htb]
  \begin{center}
    \begin{tabular}{cc}
      $T=0$ & $T=0.1$\\
\includegraphics[height=4.5cm]{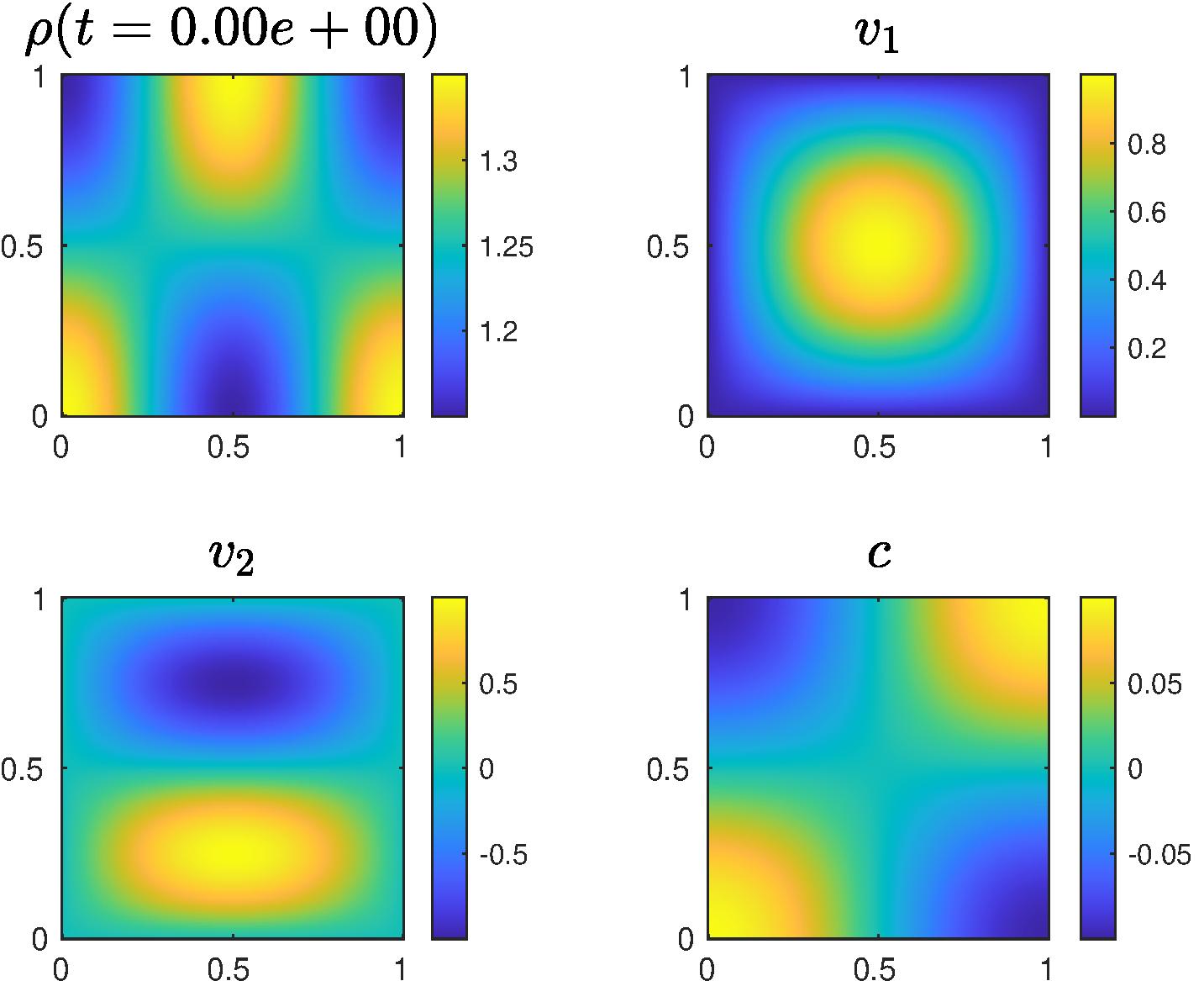}&
                                                                  \includegraphics[height=4.5cm]{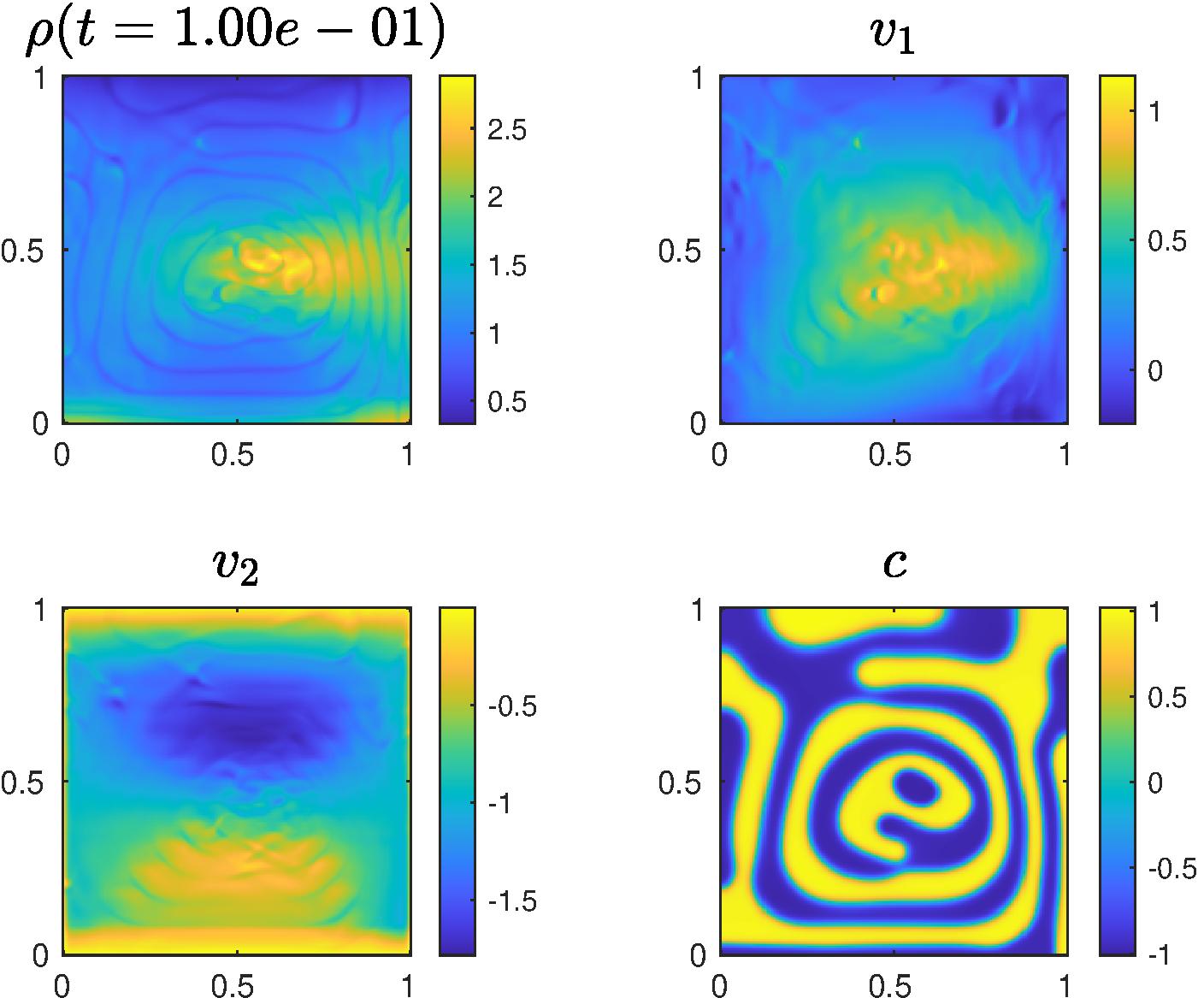}
    \end{tabular}                                                                  
  \end{center}
  \caption{Results for Test 1. Left: Initial condition, with
    $c$-variable inside spinodal region; Right: Results for $T=0.1$,
    where density increases at the bottom and separation is clearly
    visible in the $c$-variable. }
  \label{fig:test1a}
\end{figure}

\begin{figure}[htb]
  \begin{center}
    \begin{tabular}{cc}
      $T=0.3$ & $T=0.5$\\
\includegraphics[height=4.5cm]{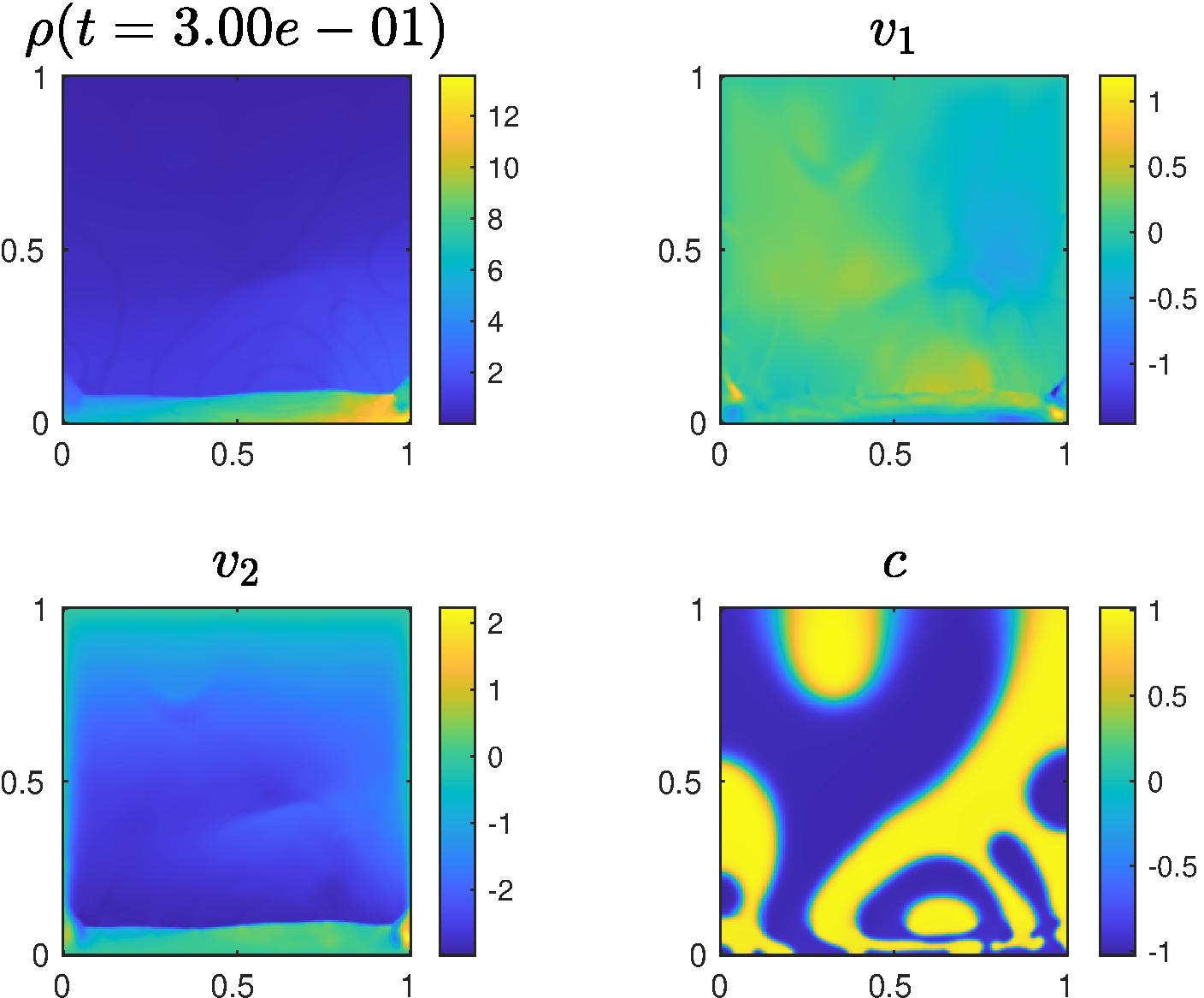}&
\includegraphics[height=4.5cm]{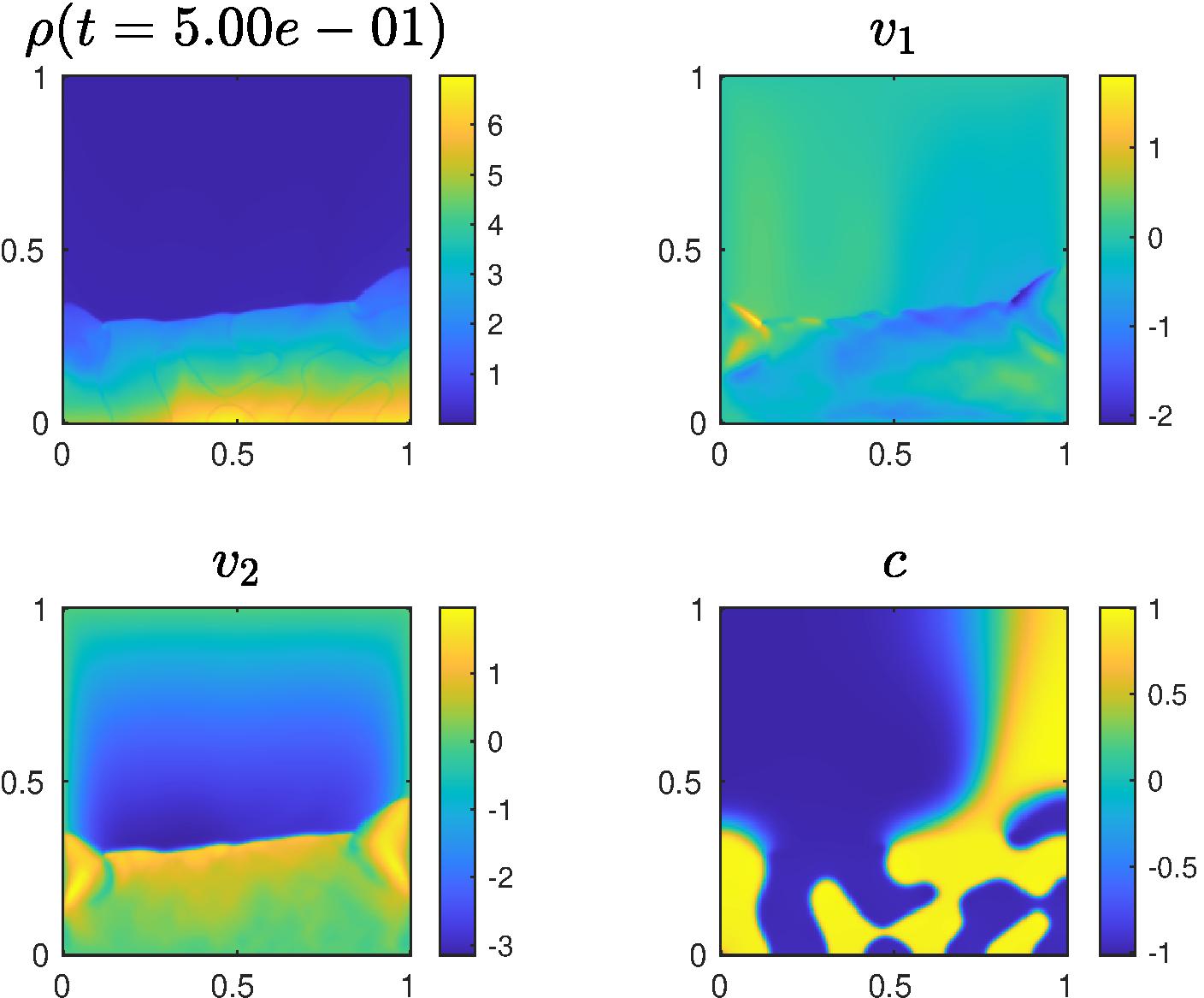}
    \end{tabular}                                                                  
  \end{center}
  \caption{Results for Test 1, $T=0.3$ (left) and $T=0.5$
    (right)     where it can be seen that density continues increasing
    at the bottom, forming an upgoing front, and nucleation     is beginning, as seen  in the $c$-variable. }
  \label{fig:test1b}
\end{figure}

\begin{figure}[htb]
  \begin{center}
    \begin{tabular}{cc}
      $T=0.7$ & $T=1.0$\\
\includegraphics[height=4.5cm]{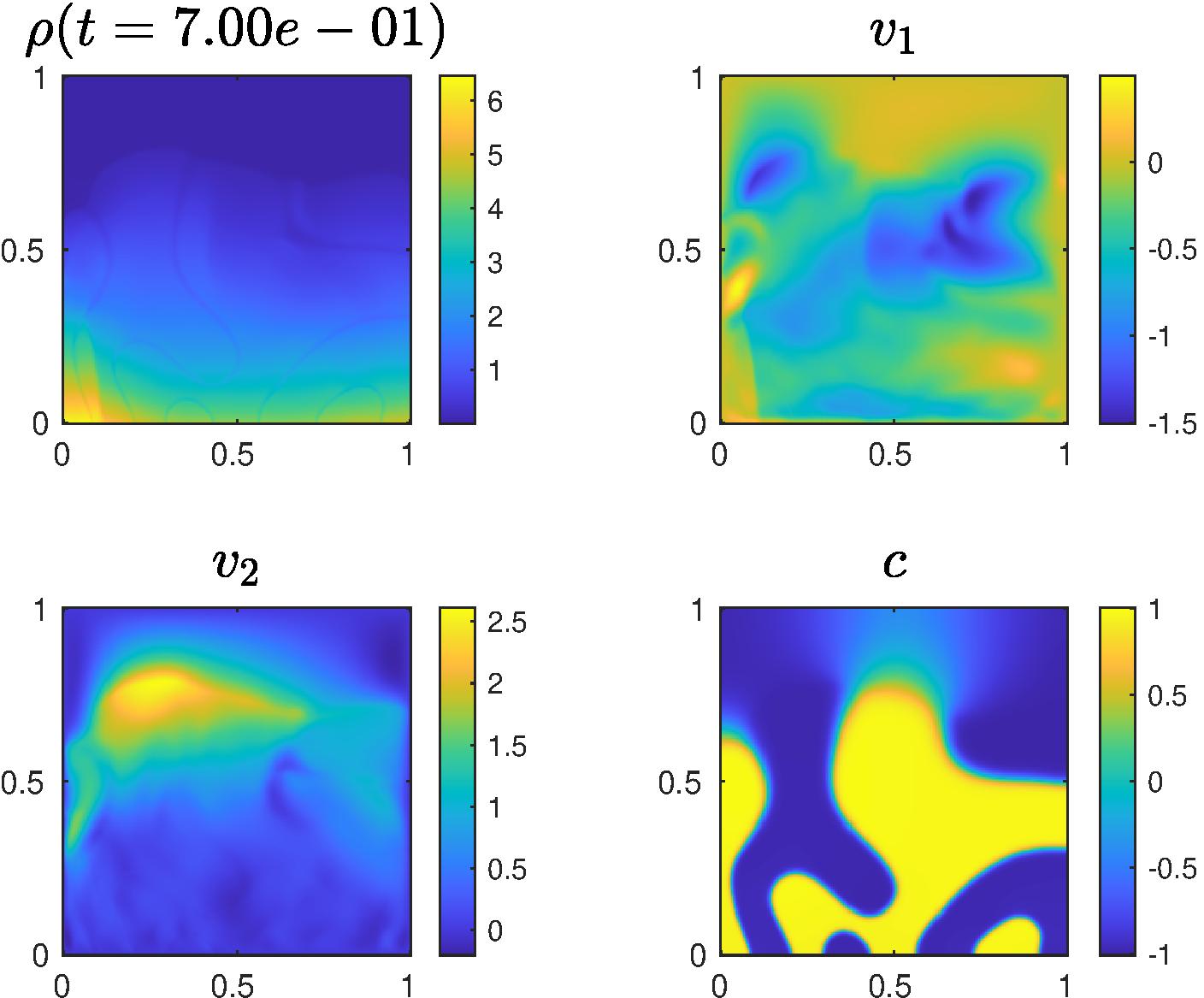}&
\includegraphics[height=4.5cm]{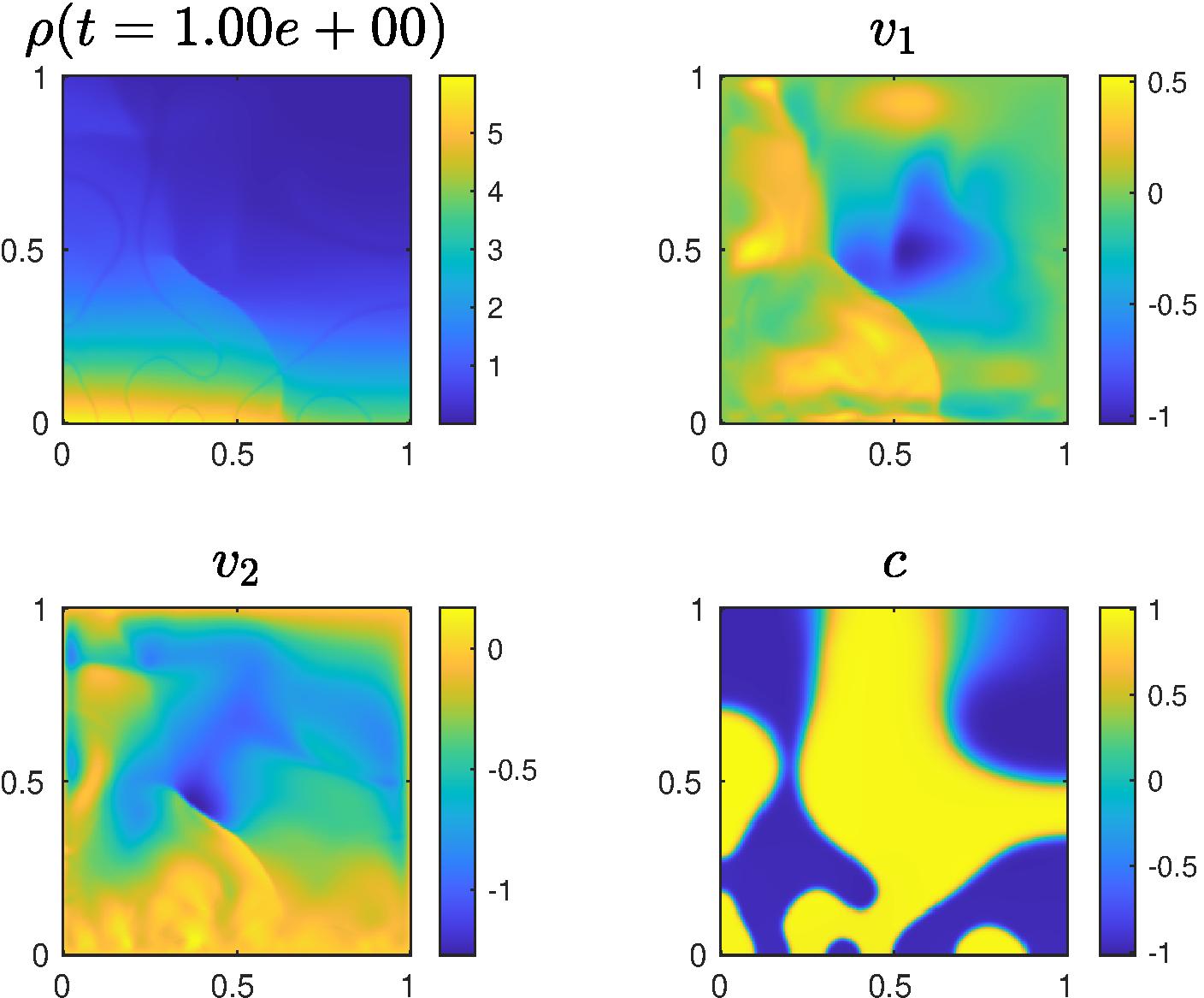}
    \end{tabular}                                                                  
  \end{center}
  \caption{Results for Test1, $T=0.7$ (left) and $T=1.0$
    (right)     where it can be seen in the velocity that vorticity
    has developed and nucleation     is increasing, as seen  in the $c$-variable. }
  \label{fig:test1c}
\end{figure}

\begin{figure}[htb]
  \begin{center}
    \begin{tabular}{cc}
      $T=0$ & $T=0.3$\\
\includegraphics[height=4.5cm]{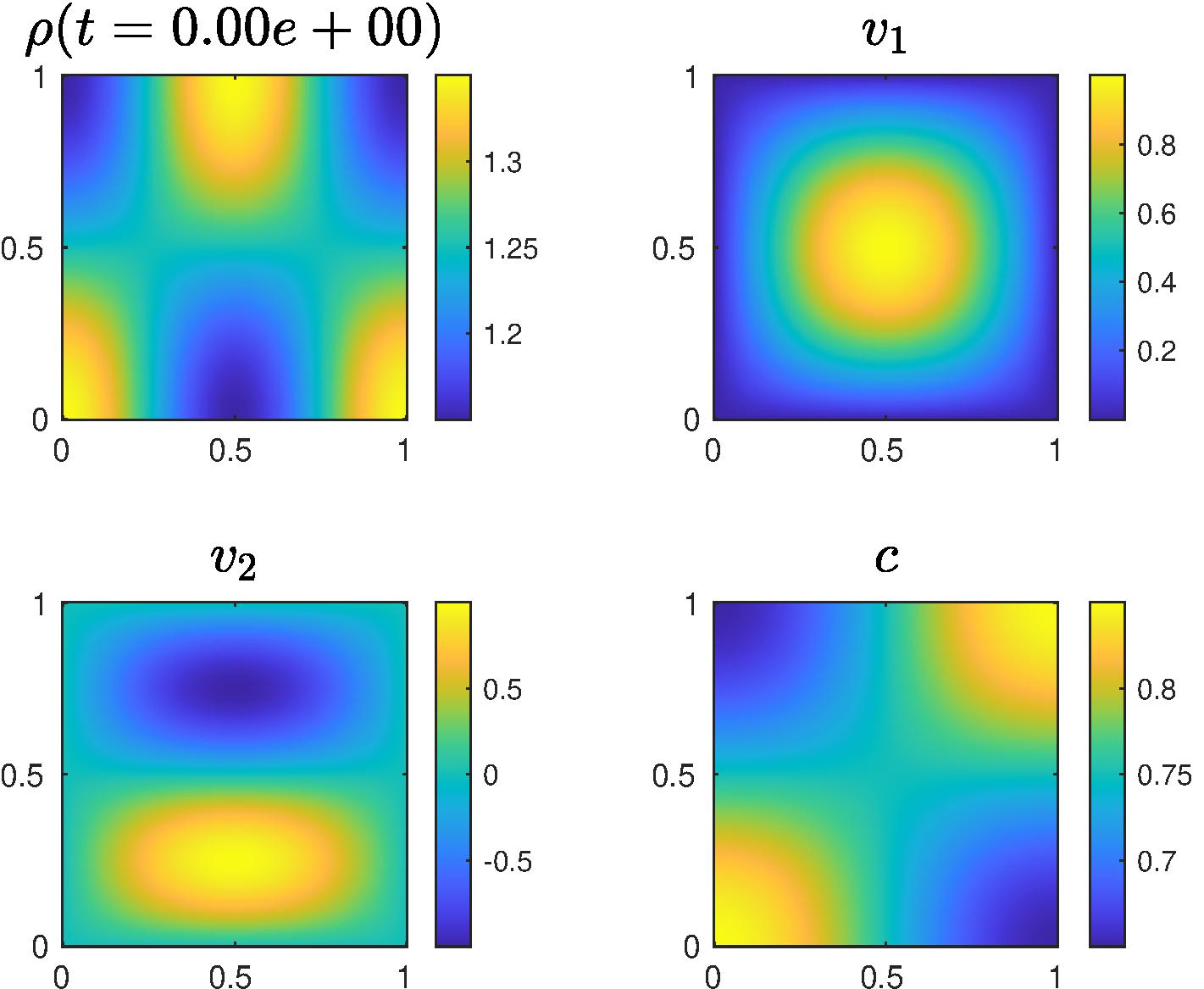}&
\includegraphics[height=4.5cm]{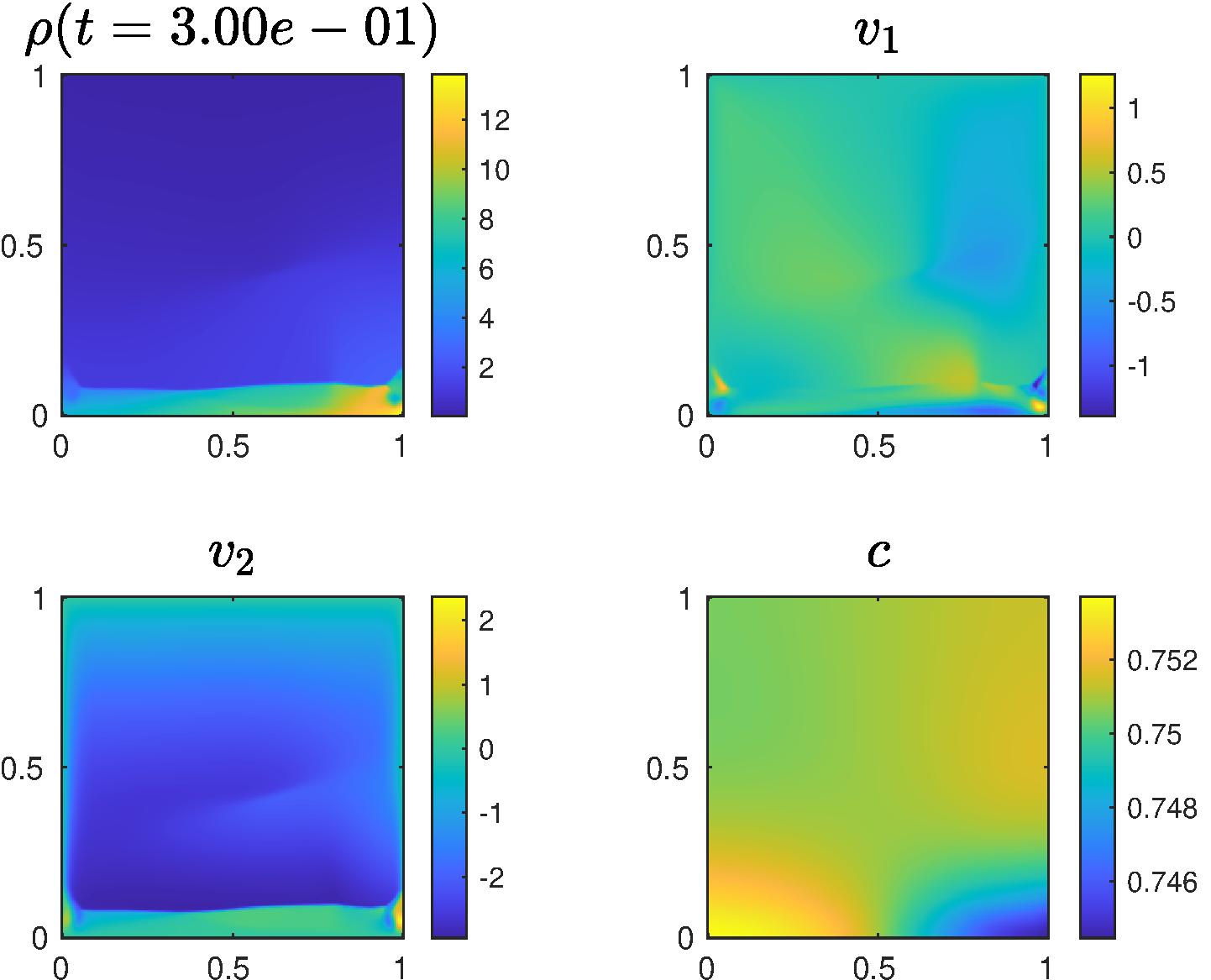}
    \end{tabular}                                                                  
  \end{center}
  \caption{Results for Test 2. Left: Initial condition, with
    $c$-variable above the spinodal region; Right: Results for $T=0.3$,
    where it can be seen that density continues increasing
    at the bottom, forming an upgoing front, while the $c$-variable is
  converging towards 0.75.}
  \label{fig:test2a}
\end{figure}

\begin{figure}[htb]
  \begin{center}
    \begin{tabular}{cc}
      $T=0.6$ & $T=1.0$\\
\includegraphics[height=4.5cm]{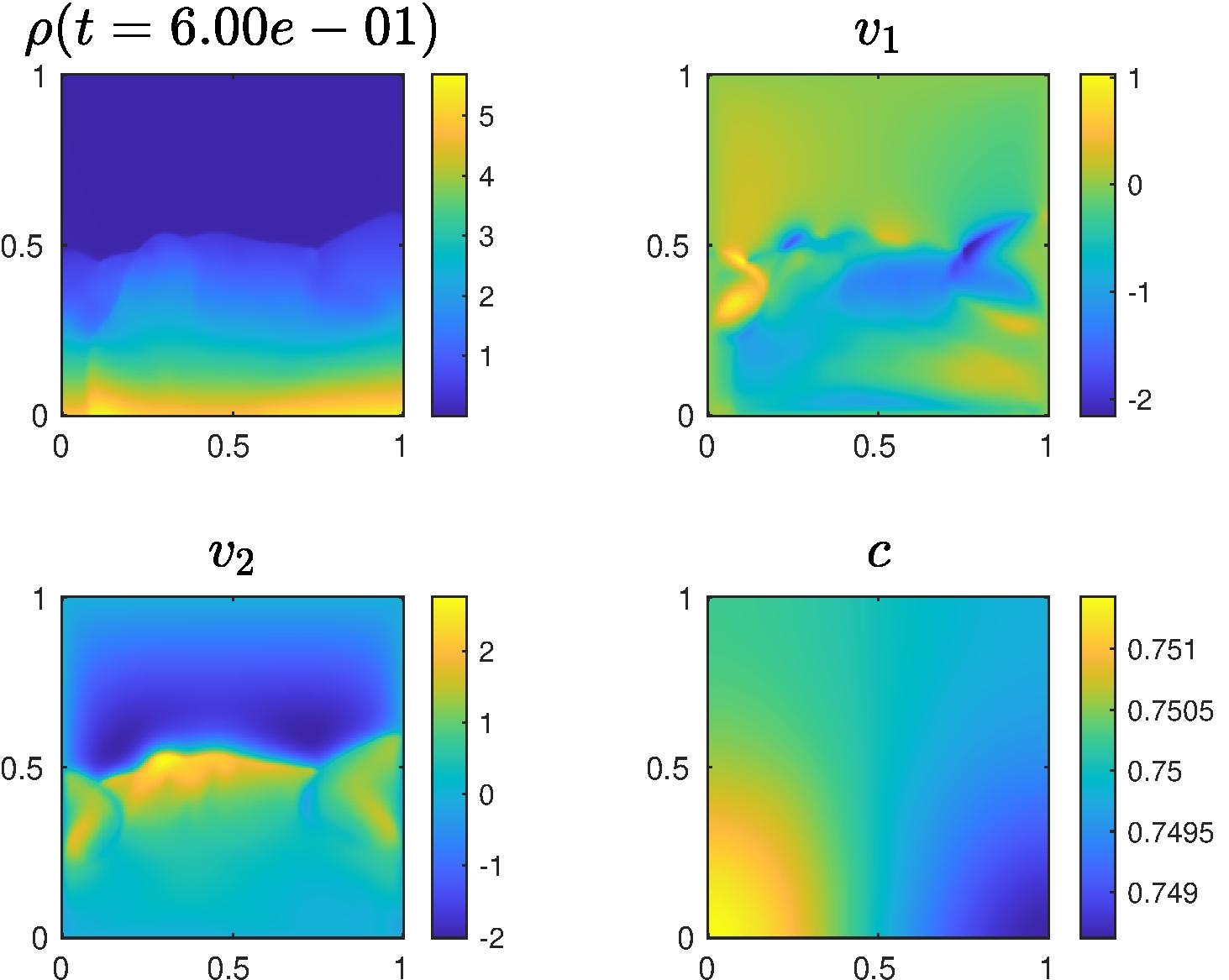}&
                                                                  \includegraphics[height=4.5cm]{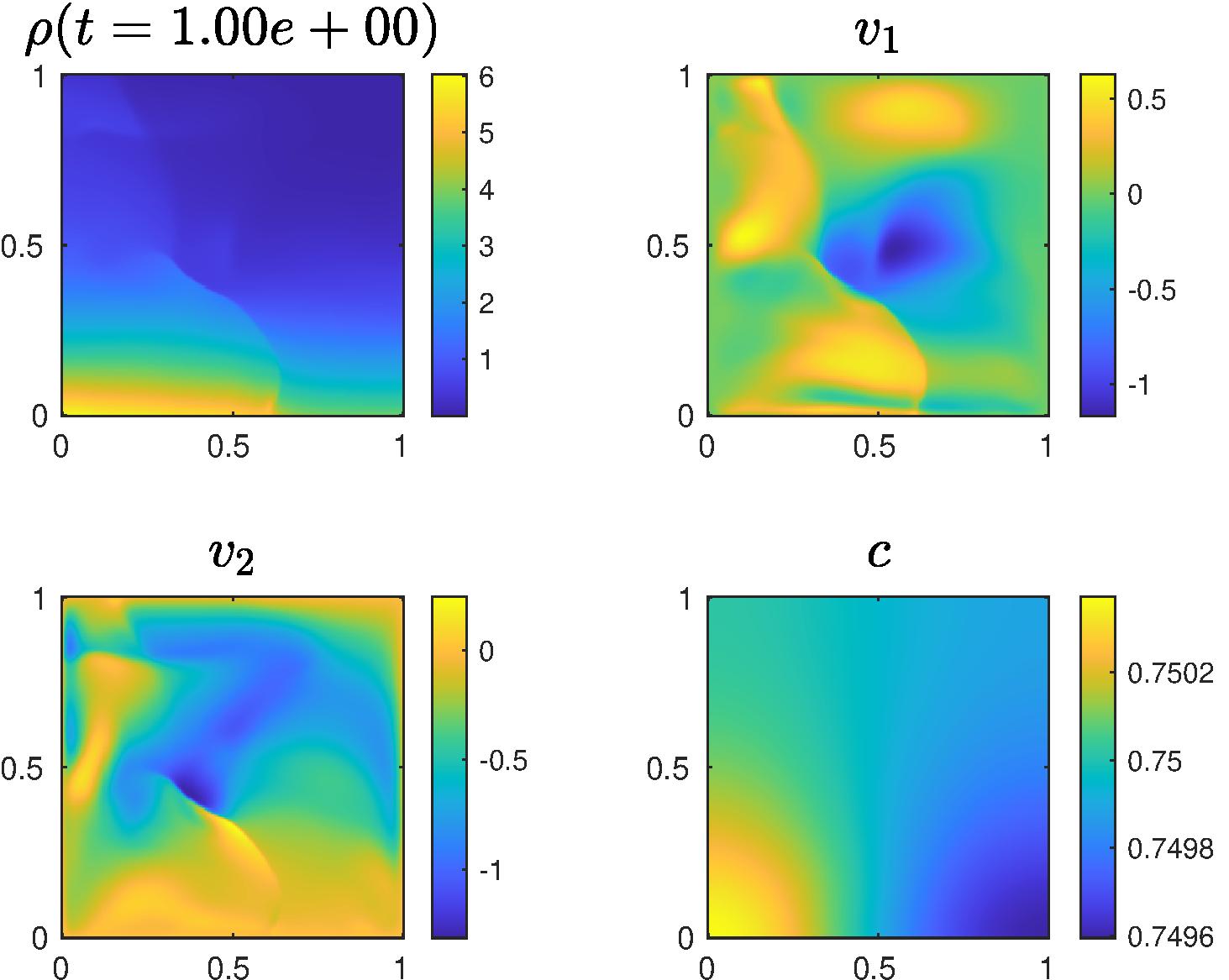}
    \end{tabular}                                                                  
  \end{center}
  \caption{Results for Test 2, $T=0.6$ (left) and $T=1.0$ (right)
    where it can be seen  in the velocity that vorticity
    has developed, while the $c$-variable has almost fully converged to 0.75.}
  \label{fig:test2b}
\end{figure}

\begin{figure}[htb]
  \begin{center}
    \begin{tabular}{cc}
      $T=0$ & $T=0.01$\\
\includegraphics[height=4.5cm]{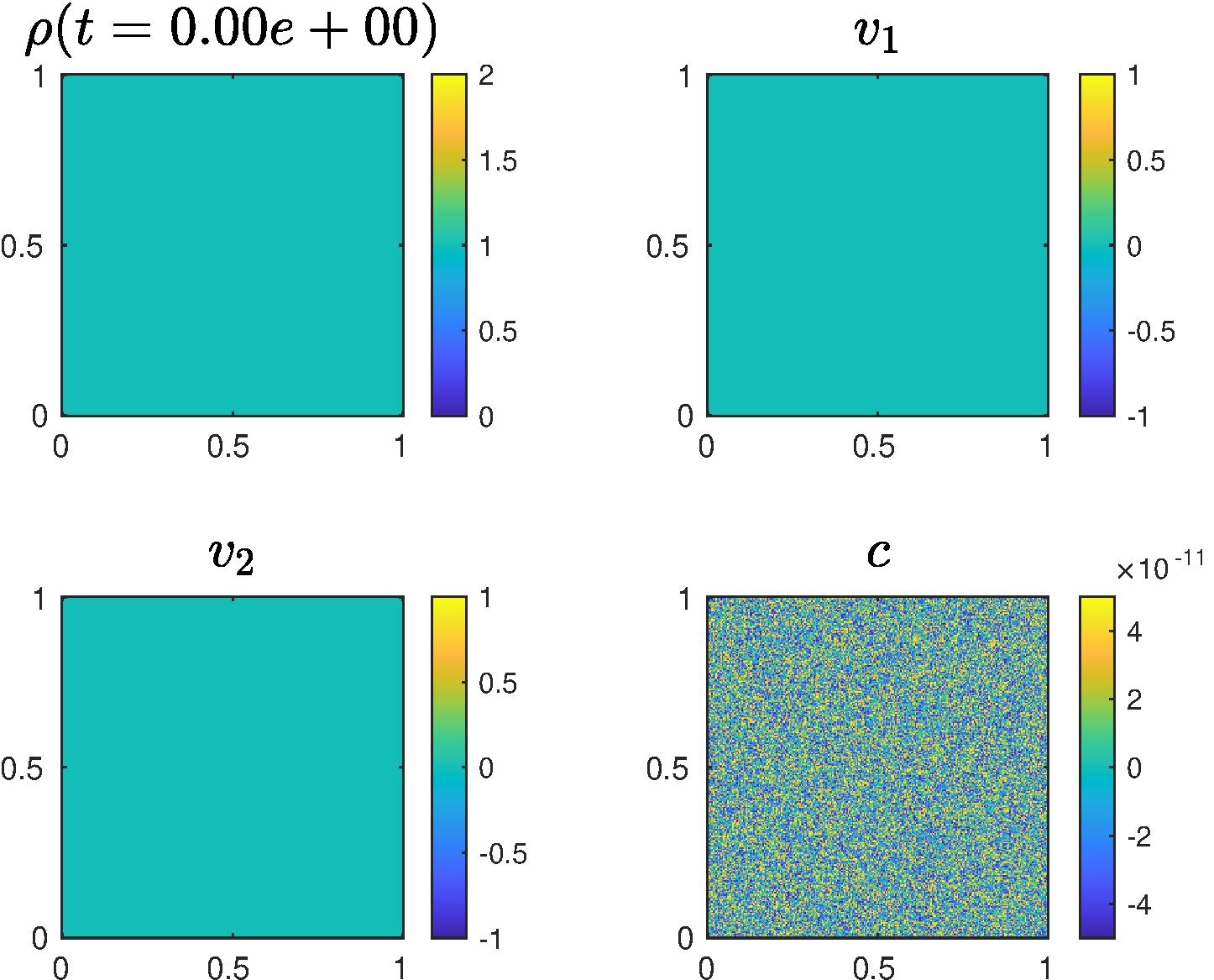}&
\includegraphics[height=4.5cm]{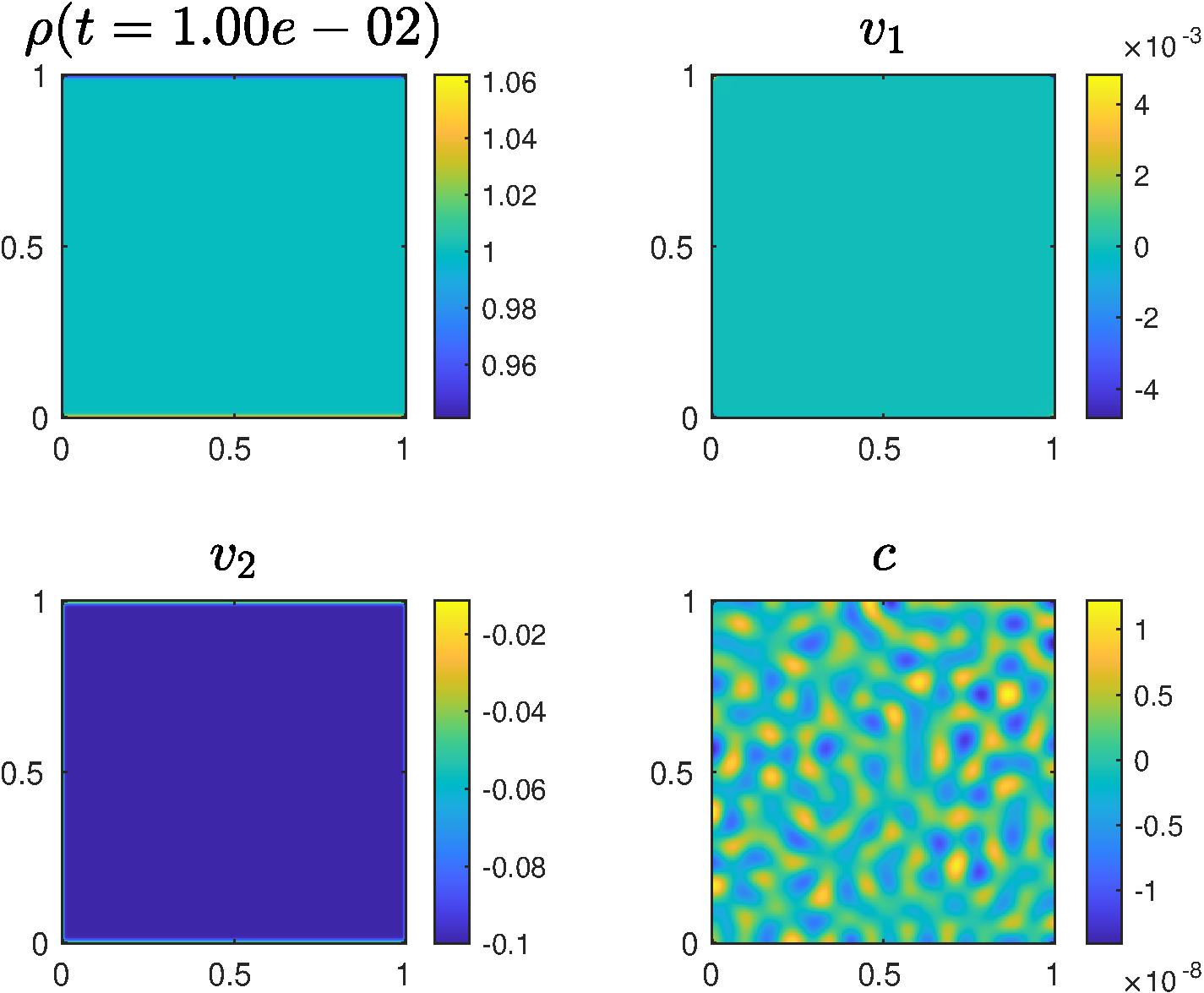}
    \end{tabular}                                                                  
  \end{center}
  \caption{Test 3. Left: Initial condition, where   it can
    be seen that the    $c$-component is in the range $(-5\, 10^{-11},
    t\, 10^{-11})$, thus lies entirely
within the spinodal region $(-\frac{1}{\sqrt3},
\frac{1}{\sqrt3})$; Right: for  $T=0.01$ onset of a  spinodal
decomposition begins appearing with an amplitude around $10^{-8}$.}
  \label{fig:test3a}
\end{figure}

\begin{figure}[htb]
  \begin{center}
    \begin{tabular}{cc}
      $T=0.02$ & $T=0.03$\\
\includegraphics[height=4.5cm]{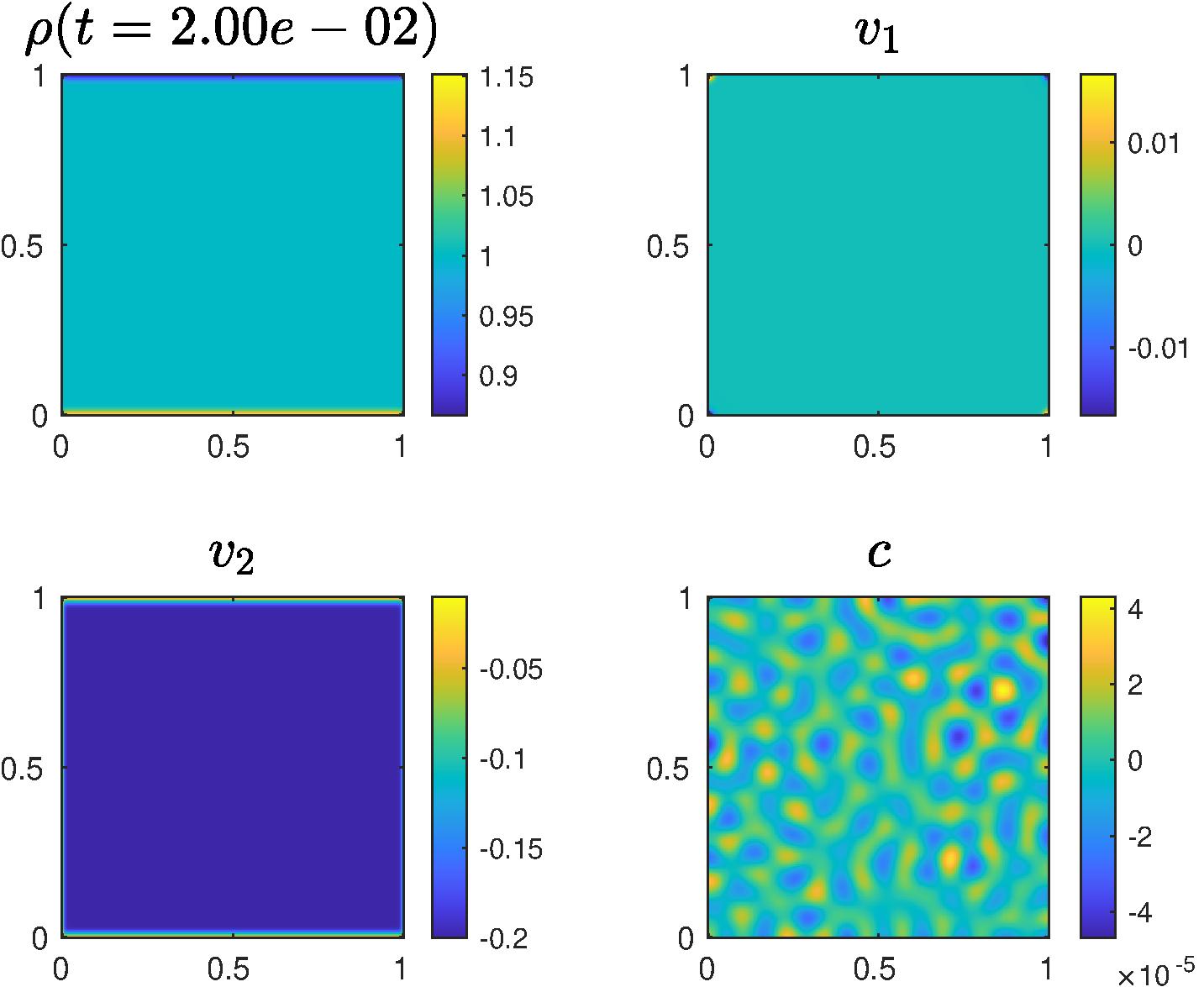}&
\includegraphics[height=4.5cm]{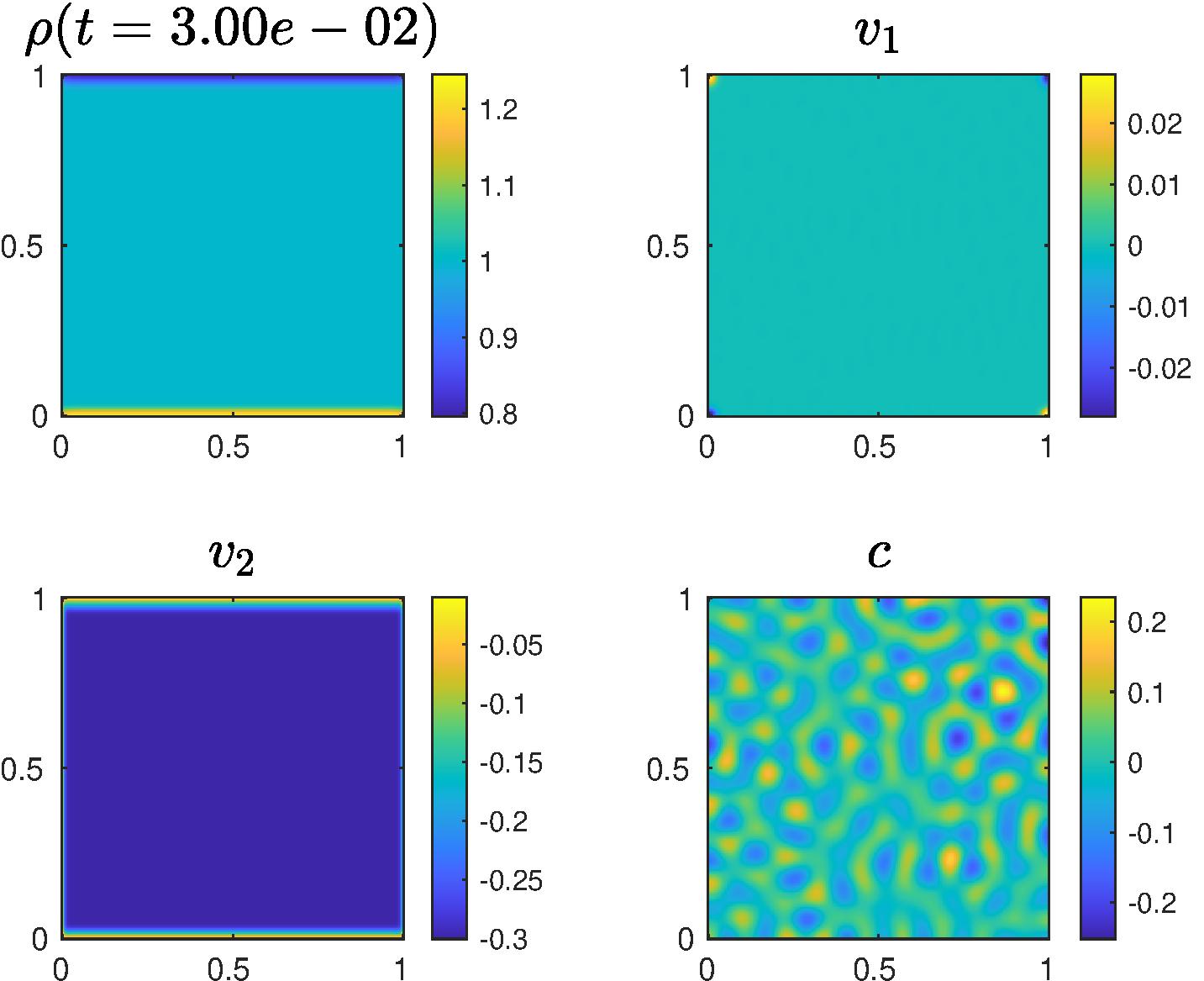}
    \end{tabular}                                                                  
  \end{center}
  \caption{Results for Test 3, where    it can 
    be seen that density increases at the bottom and the  spinodal
decomposition continues its development with an amplitude around $4\,
10^{-5}$ for 
$T=0.02$ (left) and $0.2$ for $T=0.03$ (right).}
  \label{fig:test3b}
\end{figure}

\begin{figure}[htb]
  \begin{center}
    \begin{tabular}{cc}
      $T=0.04$ & $T=0.14$\\
\includegraphics[height=4.5cm]{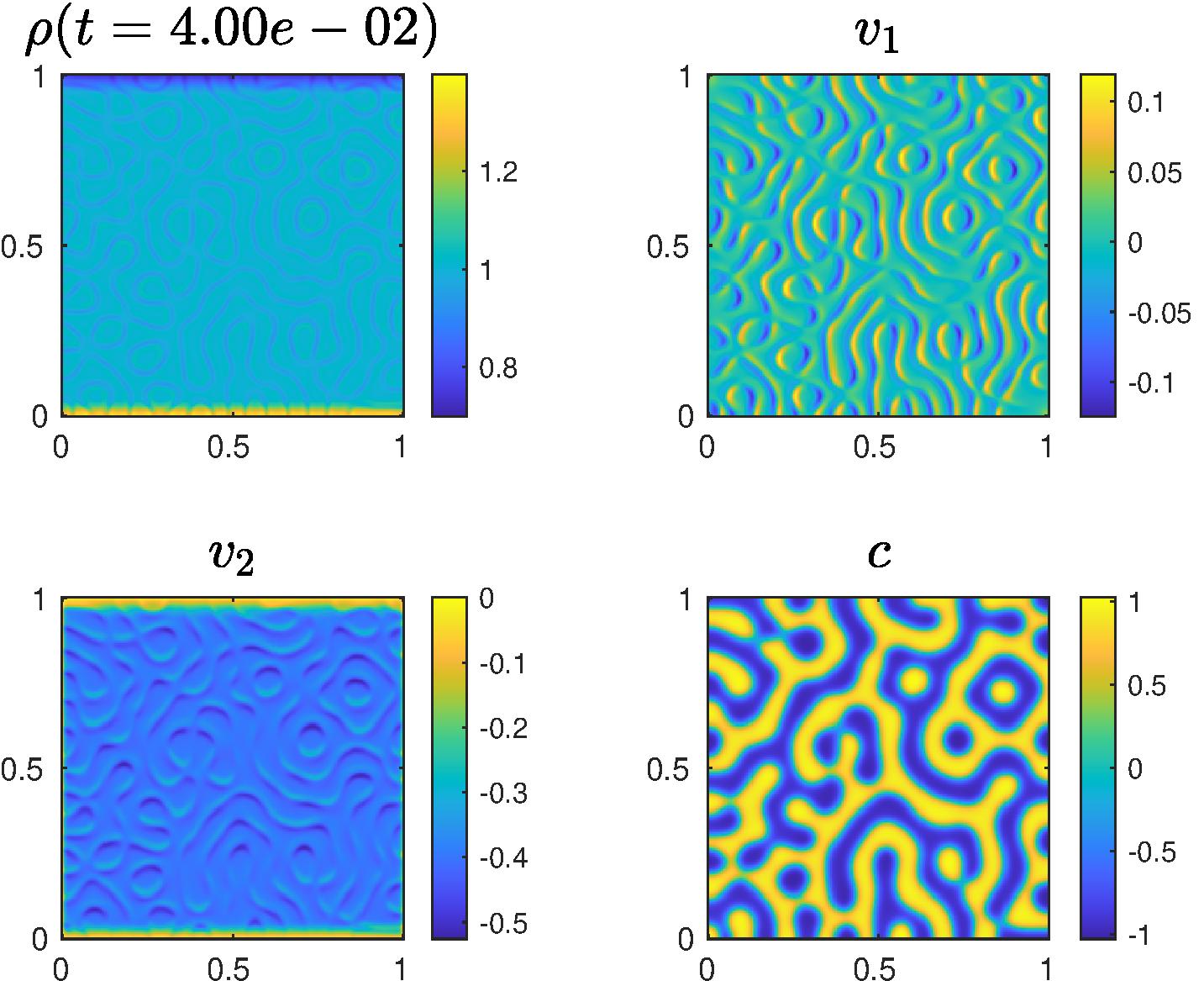}&
\includegraphics[height=4.5cm]{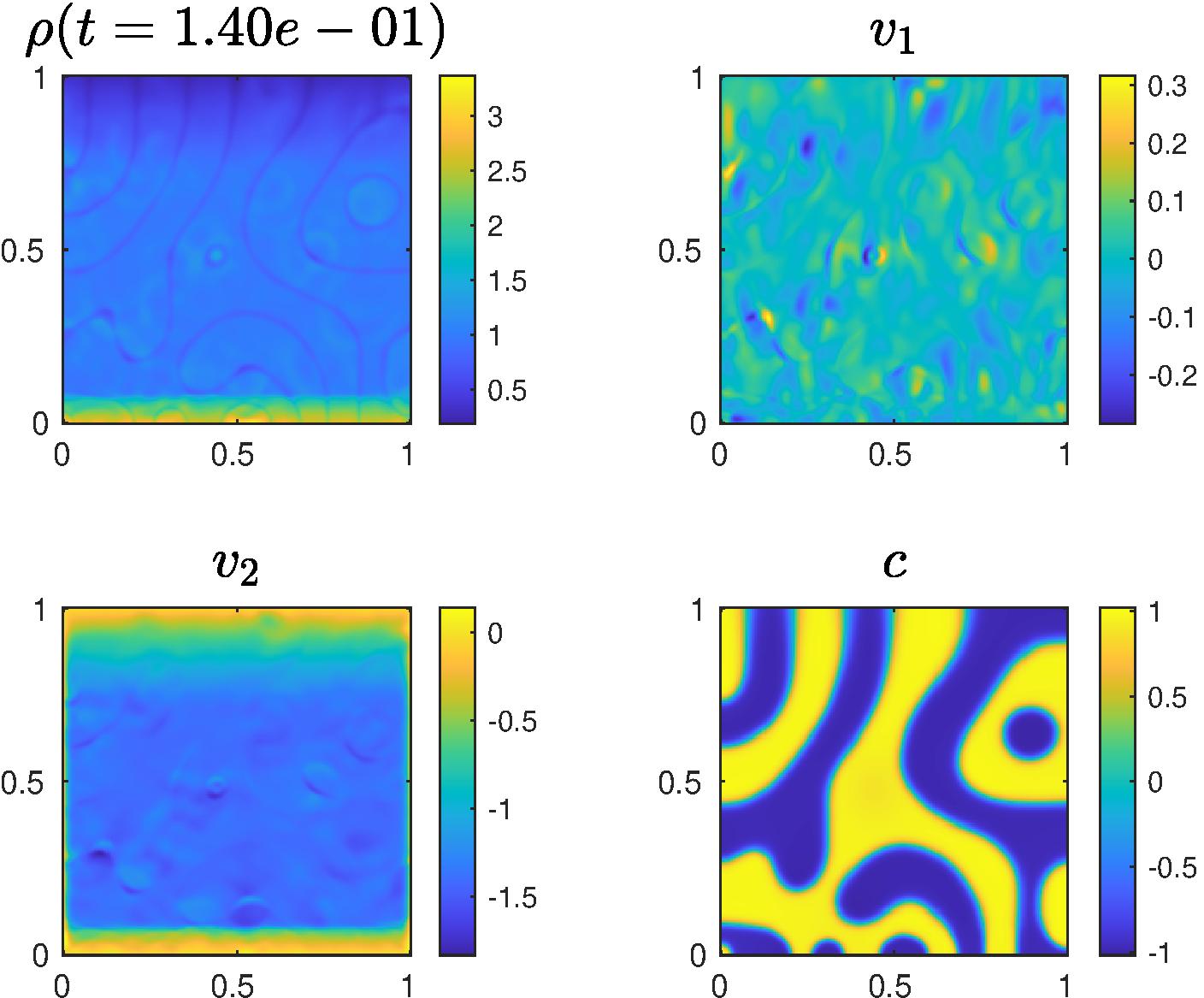}
    \end{tabular}                                                                  
  \end{center}
  \caption{Results for Test 3, where    it can 
    be seen that density continues increasing at the bottom and the  spinodal
decomposition has given a fully developed separation pattern for
$T=0.04$ (left), which has continued towards the onset of nucleation
for $T=0.14$ (right).}
  \label{fig:test3c}
\end{figure}

\begin{figure}[htb]
  \begin{center}
    \begin{tabular}{cc}
      $T=0.23$ & $T=0.24$\\
\includegraphics[height=4.5cm]{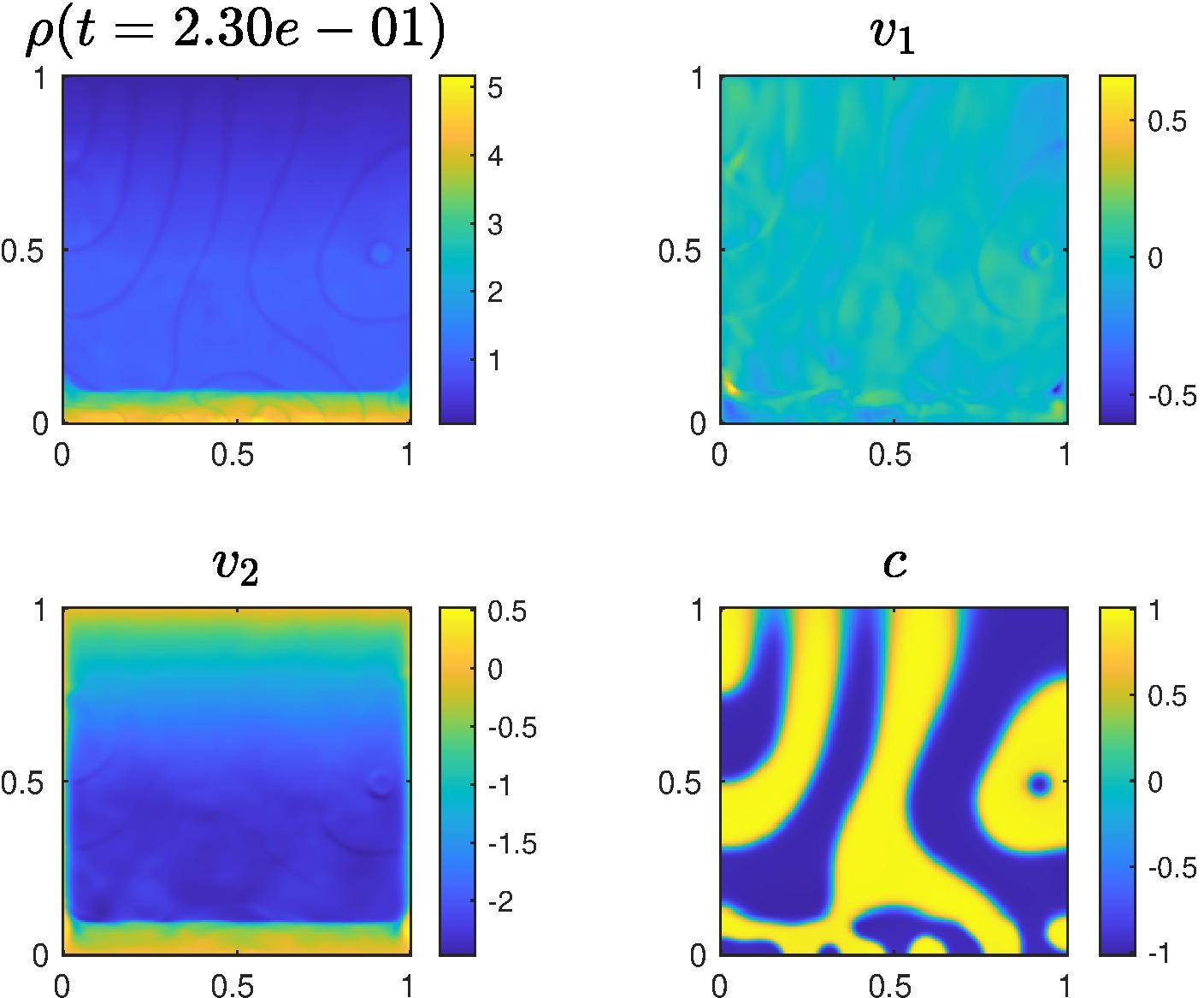}&
\includegraphics[height=4.5cm]{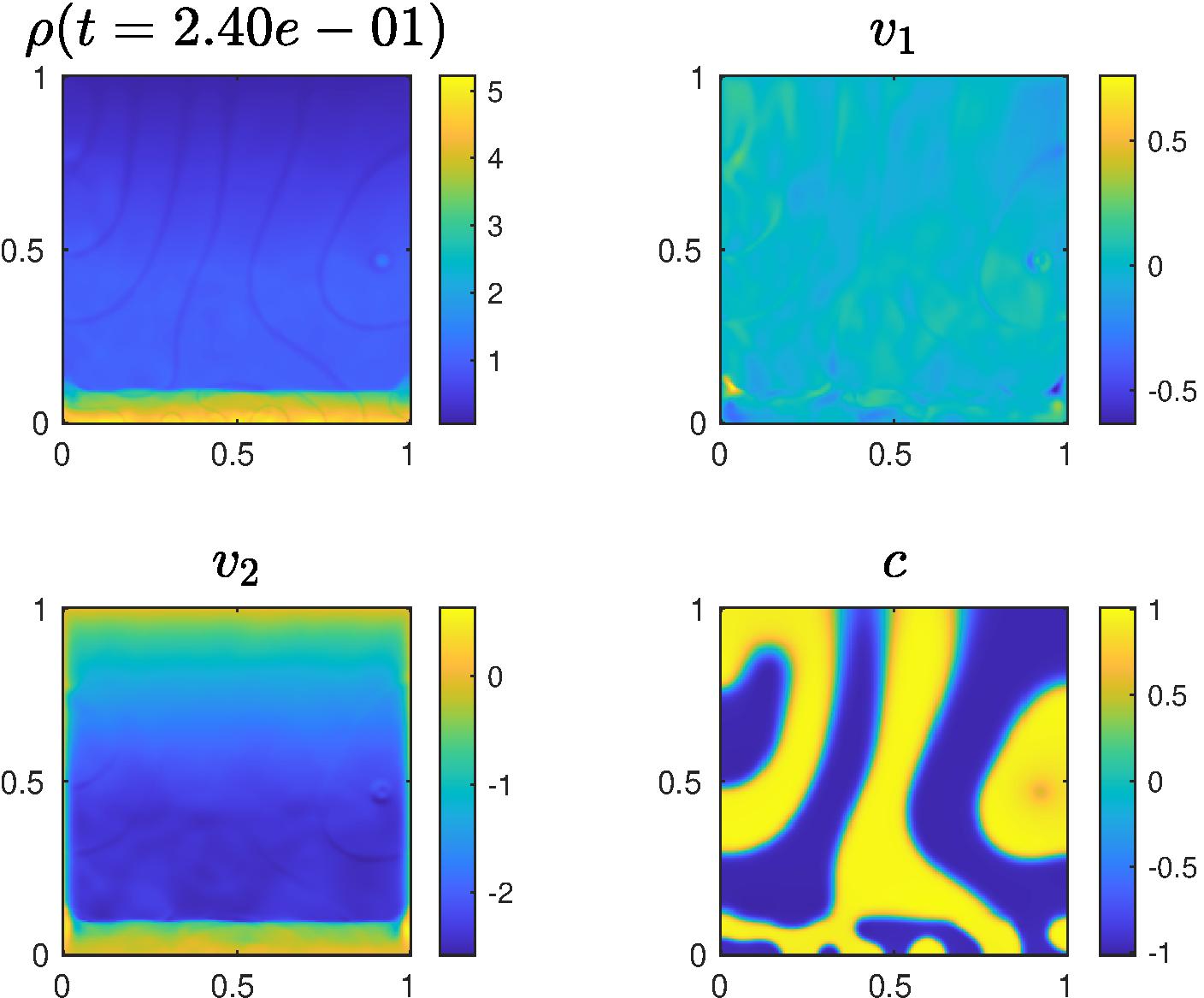}
    \end{tabular}

    \begin{tabular}{cc}
      $T=0.28$ & $T=0.29$\\
\includegraphics[height=4.5cm]{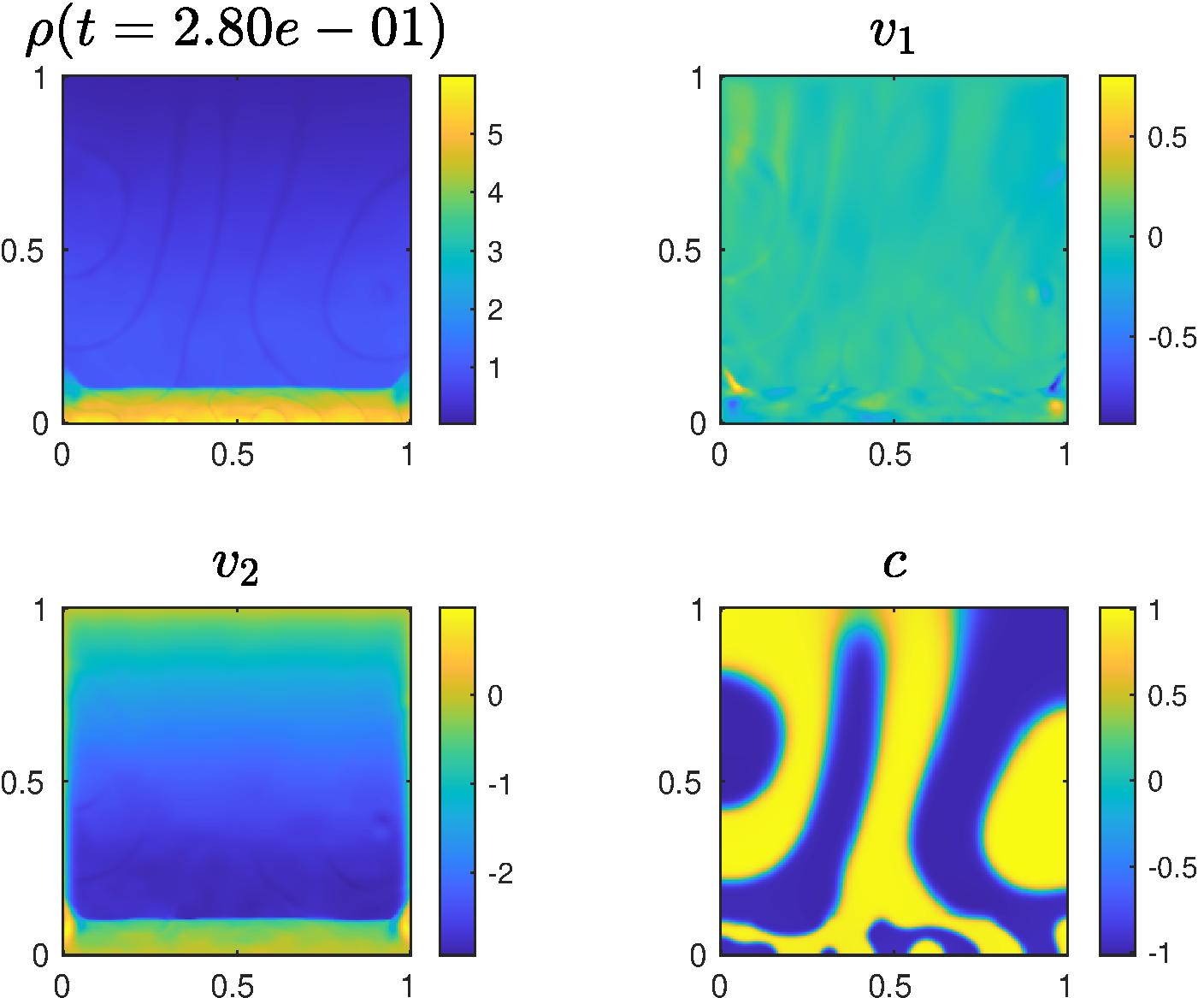}&
\includegraphics[height=4.5cm]{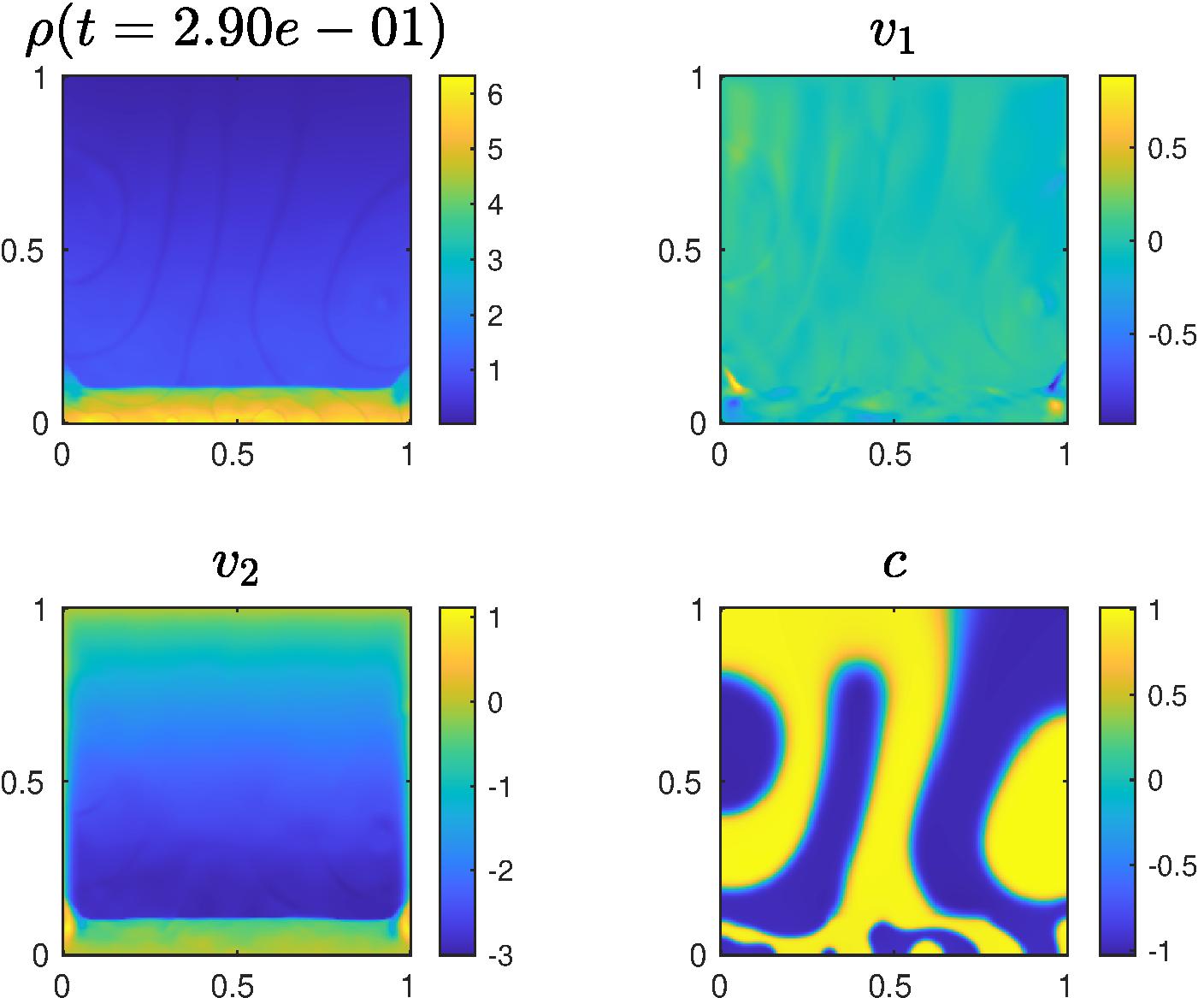}
    \end{tabular}                                                                  
  \end{center}
  \caption{Results for Test 3, where  nucleation is visible in the
    $c$-variable near $(0.9, 0.5)$, with an ever-shrinking region of
    the second fluid disappearing.}
  \label{fig:test3e}
\end{figure}

We show in Figures \ref{fig:test1a}-\ref{fig:test3e} some snapshots of
the results obtained for all the tests with *-DIRKSA, $M=256$, $G=-10$, 
$\nu=10^{-3}, \lambda=10^{-4}, \varepsilon=10^{-4}$, which correspond to 
flows with Reynolds number roughly in the range $[10^2, 10^3]$.

In Figure \ref{fig:test1a}, it can be seen that 
the $c$-component in the initial condition for Test 1 lies entirely
within the spinodal region $(-\frac{1}{\sqrt3},
\frac{1}{\sqrt3})$. Therefore, at the early stages of the 
simulation, at $T=0.1$, separation occurs forming complex
patterns. Meanwhile, it can be appreciated from the pictures
corresponding to the $\rho$-variable that  gravity is acting so that
density increases at the bottom boundary, $y=0$. This trend continues
in Figure \ref{fig:test1b}, where it can be appreciated that the
maximal density
continues increasing, a sharp
upgoing front develops for all the variables, but $c$, where the
diffuse interface experiments many topological changes, with growing
regions, a phenomenon named nucleation. At the final stages of the
simulation, as seen
in Figure \ref{fig:test1c}, the maximal density decreases, nucleation
continues and the bulk flow  enters into a seemingly turbulent regime.

In Figure \ref{fig:test2a}, it can be seen that 
the $c$-component in the initial condition for Test 2 lies above
 the spinodal region. Therefore, as seen in the pictures for the
 $c$-variable in the snapshots shown in this Figure and also in Figure
 \ref{fig:test2b}, the fluid remains almost homogeneous, with $c$
 tending to $3/4$ in the whole domain, which corresponds to  mass
 fractions $c_1=\frac{7}{8}, c_{2}=\frac{1}{8}$, which are exactly
 the initial proportions of the individual species. This fact means
 that the rest of the equations behave like a uniform fluid governed
 by the compressible Navier-Stokes equations under gravitation, with
 Reynolds number high enough for a seemingly turbulent regime.

In Figure \ref{fig:test3a}, it can be seen that the  
 $c$-component in the initial condition for Test 3 is almost 0, thus lies entirely
within the spinodal region $(-\frac{1}{\sqrt3},
\frac{1}{\sqrt3})$. Therefore, at the early stages of the simulation, at $T=0.01$,
a typical spinodal decomposition begins appearing in the form of a 
medium-frequency pattern, corresponding to a solution as in
\eqref{eq:solspin0}-\eqref{eq:solspin}, for $k_1, k_2\in \mathbb{N}$ that
minimize
\begin{align*}
  -(k_1^2+k_2^2)+\varepsilon\pi^2 (k_1^2+k_2^2)^2,
\end{align*}
corresponding to the expression in \eqref{eq:solspin}, for $c_0=0$,
taking into account that $\psi''(c_0)=-1$.

The rest of the simulation can be seen in Figures
\ref{fig:test3b}-\ref{fig:test3e}, where the spinodal decomposition
continues until nucleation. Density increases at the bottom, but no
clear turbulence is appreciated, maybe due to the  simulation having
been carried only until $T=0.29$.

\subsection{Linear solvers}
 The multigrid solver mentioned in subsection
\ref{ss:solver} for the solution of \eqref{eq:227}  has shown
 a very satisfactory performance throughout all tests, requiring almost
always only one iteration to achieve convergence to double precision
digits.

The performance of the solvers  for the solution of \eqref{eq:226} is
more complex. In tables \ref{tbl:nits7}, \ref{tbl:nits6} and \ref{tbl:nits0} we show the average number of multigrid iterations to solve \eqref{eq:226} for tests 1, 2, 3, respectively, with $\nu=10^{-2},10^{-3},10^{-4}$ and $\varepsilon=10^{-3},10^{-4},10^{-5}$, for a stopping criterion based on a relative decrease of the residual by a factor of $10^{-6}$.

It can be deduced from these tables that the number of iterations grows slowly with $M$ for all cases.

In Table \ref{tbl:cmpcgmg} we show a comparison between the multigrid
 and the preconditioned conjugate gradient solvers in subsection \ref{ss:solver}.
 The results have been obtained with *-DIRKSA, $M=16,32,64,128,256$, $G=-10$,
$\nu=10^{-2}, \lambda=10^{-3}, \varepsilon=10^{-4}$. It can be deduced  that
the preconditioned conjugate gradient solver uses less CPU time than
the multigrid solver, although the latter takes fewer iterations than
the former.

\begin{table}[htb]
  \addtolength{\tabcolsep}{-2pt}
  \begin{center}
\begin{tabular}{|c|ccc|ccc|ccc|}
  \hline
  &\multicolumn{3}{c|}{$\varepsilon=10^{-3}$}
  &\multicolumn{3}{c|}{$\varepsilon=10^{-4}$}
  &\multicolumn{3}{c|}{$\varepsilon= 10^{-5}$}\\\hline
$M/\nu$  &$10^{-1}$&$10^{-2}$&$10^{-3}$&
  $10^{-1}$&$10^{-2}$&$10^{-3}$&
  $10^{-1}$&$10^{-2}$&$10^{-3}$\\  \hline
 16& 2.4 & 2.8 & 3.3  & 1.9 & 2.5 & 3.1 & 1.8 & 1.9 & 2.1\\
 32& 3.5 & 4.1 & 5.8  & 2.3 & 3.2 & 3.8 & 2.0 & 2.0 & 2.7\\
 64& 4.2 & 5.2 & 4.5  & 3.6 & 3.8 & 4.2 & 2.2 & 2.3 & 3.0\\
128& 5.1 & 6.4 & 5.6  & 4.4 & 4.5 & 5.7 & 3.1 & 3.7 & 3.8\\
256& 5.9 & 7.3 & 6.3  & 4.8 & 5.6 & 6.8 & 3.8 & 3.9 & 4.5\\
  \hline
\end{tabular}
\end{center}
\caption{Average number of multigrid iterations to solve \eqref{eq:226} for Test 1.}
\label{tbl:nits7}
\end{table}

\begin{table}[htb]
  \addtolength{\tabcolsep}{-2pt}
\begin{center}
\begin{tabular}{|c|ccc|ccc|ccc|}
  \hline
  &\multicolumn{3}{c|}{$\varepsilon=10^{-3}$}
  &\multicolumn{3}{c|}{$\varepsilon=10^{-4}$}
  &\multicolumn{3}{c|}{$\varepsilon=10^{-5}$}\\\hline
$M/\nu$  &$10^{-1}$&$10^{-2}$&$10^{-3}$&
  $10^{-1}$&$10^{-2}$&$10^{-3}$&
  $10^{-1}$&$10^{-2}$&$10^{-3}$\\  \hline
16 & 1.2 & 1.2 & 1.3 & 1.1 & 1.1 & 1.2 & 1.1 & 1.1 & 1.1\\
 32& 1.3 & 1.3 & 1.6 & 1.1 & 1.1 & 1.3 & 1.1 & 1.1 & 1.2\\
 64& 1.4 & 1.6 & 2.2 & 1.1 & 1.2 & 1.3 & 1.1 & 1.3 & 1.4\\
128& 1.6 & 1.9 & 2.9 & 1.2 & 1.2 & 1.2 & 1.2 & 1.4 & 1.4\\
256& 1.8 & 2.3 & 3.4 & 1.3 & 1.3 & 1.3 & 1.3 & 1.4 & 1.4  
  \\
  \hline
\end{tabular}
\end{center}
\caption{Average number of multigrid iterations to solve \eqref{eq:226} for Test 2.}
\label{tbl:nits6}
\end{table}

\begin{table}[htb]
  \addtolength{\tabcolsep}{-2pt}
  \begin{center}
    \begin{tabular}{|c|ccc|ccc|ccc|}
  \hline
   &\multicolumn{3}{c|}{$\varepsilon=10^{-3}$}
  &\multicolumn{3}{c|}{$\varepsilon=10^{-4}$}
  &\multicolumn{3}{c|}{$\varepsilon= 10^{-5}$}\\\hline
$M/\nu$ &$10^{-1}$&$10^{-2}$&$10^{-3}$&
  $10^{-1}$&$10^{-2}$&$10^{-3}$&
  $10^{-1}$&$10^{-2}$&$10^{-3}$\\  \hline
 16& 3.3 & 3.1 & 3.1 & 2.6 & 2.3 & 2.4 & 2.6 & 2.3 & 2.3\\
 32& 4.5 & 4.4 & 4.5 & 3.3 & 2.9 & 3.0 & 3.1 & 2.7 & 2.8\\
 64& 6.4 & 5.6 & 5.5 & 4.5 & 3.7 & 3.7 & 4.3 & 3.8 & 3.4\\
128& 6.5 & 6.0 & 5.5 & 6.4 & 5.7 & 5.3 & 5.4 & 4.4 & 4.0\\
256& 7.0 & 6.0 & 5.5 & 7.0 & 6.0 & 5.5 & 7.2 & 5.7 & 5.4 
   \\
  \hline
\end{tabular}
\end{center}
\caption{Average number of multigrid iterations to solve \eqref{eq:226} for Test 3.}
\label{tbl:nits0}
\end{table}

\begin{table}[htb]
  \begin{center}
  \addtolength{\tabcolsep}{-6pt}
\begin{tabular}{|c|cc|cc|}
  \hline
  &\multicolumn{2}{c|}{CG}
  &\multicolumn{2}{c|}{MG}\\
  \hline
  $M$&avg. its.&CPU&avg. its.&CPU\\\hline
16 &10.00&0.05&    2.96&0.07\\  
32 &10.48&0.30&    3.00&0.40\\  
64 &12.27&2.22&    3.20&2.53\\  
128&15.61&15.46&   4.54&19.53\\ 
  256&18.05&137.27&  5.96&195.80\\
  \hline
\end{tabular}
\end{center}
\caption{Comparison of average number of iterations and CPU time for
  the preconditioned conjugate gradient and multigrid solvers for \eqref{eq:226} and Test 1.}
\label{tbl:cmpcgmg}
\end{table}

\section{Conclusions and future work}\label{s:conclusion}
In this paper, we propose efficient linearly implicit-explicit schemes for the
two-dimensional compressible isentropic Cahn-Hilliard-Navier-Stokes equations.
Some tests are performed to show that they  achieve second-order
accuracy under time-step stability restrictions dictated only by the
convective part of the equations.

As future research, we plan to extend these techniques to other, 
stiffer pressure laws, and to  a
three-dimensional setting with Galerkin techniques. We also plan the
extension  of these techniques to quasi-incompressible models (see
\cite{LT98}).

A crucial part of the algorithms is the iterative linear solvers used
for solving the system related to the Cahn-Hilliard subequation. We
plan to explore the possibility of using the multigrid solver as
preconditioner for the conjugate gradient solver.

\section*{Acknowledgments}
I wish to express my gratitude to Raimund Bürger, from the  University
of Concepción, Chile, for suggesting to look at Siano's
paper \cite{Siano79} and to Rafael Ordóñez for preliminary work
on Cahn-Hilliard simulations.

This paper has received financial support from the research projects
PID2020-117211GB-I00, granted by MCIN/ AEI /10.13039/501100011033, and 
CIAICO/2021/227, granted by GVA.


\begin{thebibliography}{10}

\bibitem{AbelsFeireisl08}
Helmut Abels and Eduard Feireisl.
\newblock {On a diffuse interface model for a two-phase flow of compressible
  viscous fluids}.
\newblock {\em Indiana Univ. Math. J.}, 57(2):659--698, 2008.

\bibitem{BBMZ19a}
Antonio Baeza, Raimund Burger, Pep Mulet, and David Zorio.
\newblock {On the Efficient Computation of Smoothness Indicators for a Class of
  WENO Reconstructions}.
\newblock {\em Journal of Scientific Computing}, 80(2):1240--1263, AUG 2019.

\bibitem{BBMZ19b}
Antonio Baeza, Raimund Burger, Pep Mulet, and David Zorio.
\newblock {WENO Reconstructions of Unconditionally Optimal High Order}.
\newblock {\em SIAM Journal on Numerical Analysis}, 57(6):2760--2784, 2019.

\bibitem{BBMRV15}
Sebastiano Boscarino, Raimund B\"urger, Pep Mulet, Giovanni Russo, and Luis~M.
  Villada.
\newblock {Linearly Implicit Imex Runge-Kutta Methods for a Class of Degenerate
  Convection-Diffusion Problems}.
\newblock {\em {SIAM J. Sci. Comp.}}, {37}({2}):{B305--B331}, {2015}.

\bibitem{BrandtLivne2011}
Achi Brandt and Oren~E. Livne.
\newblock {\em {Multigrid techniques---1984 guide with applications to fluid
  dynamics}}, volume~67 of {\em Classics in Applied Mathematics}.
\newblock Society for Industrial and Applied Mathematics (SIAM), Philadelphia,
  PA, revised edition, 2011.

\bibitem{CahnHilliard59}
J.W. Cahn and J.E. Hilliard.
\newblock {Free energy of a nonuniform system .3. Nucleation in a 2-component
  incompressible fluid}.
\newblock {\em J. Chem. Phys.}, 31(3):688--699, 1959.

\bibitem{Elliott89}
C.~M. Elliott.
\newblock {The {C}ahn-{H}illiard model for the kinetics of phase separation}.
\newblock In {\em Mathematical models for phase change problems (\'{O}bidos,
  1988)}, volume~88 of {\em Internat. Ser. Numer. Math.}, pages 35--73.
  Birkh\"{a}user, Basel, 1989.

\bibitem{ElliottFrench87}
Charles~M. Elliott and Donald~A. French.
\newblock {Numerical Studies of the Cahn-Hilliard Equation for Phase
  Separation}.
\newblock {\em IMA Journal of Applied Mathematics}, 38(2):97--128, 05 1987.

\bibitem{HeShi20}
Qiaolin He and Xiaoding Shi.
\newblock {Numerical Study of Compressible Navier-Stokes-Cahn-Hilliard System}.
\newblock {\em Comm. Math. Sci.}, 18(2):571--591, 2020.

\bibitem{Jacmin99}
D.~Jacqmin.
\newblock {Calculation of two-phase Navier-Stokes flows using phase-field
  modeling}.
\newblock {\em Journal of Computational Physics}, 155(1):96--127, OCT 10 1999.

\bibitem{Kynch52}
GJ~Kynch.
\newblock {A Theory of Sedimentation}.
\newblock {\em Trans. Faraday Soc.}, 48(2):166--176, 1952.

\bibitem{LT98}
J.~Lowengrub and L.~Truskinovsky.
\newblock {Quasi-incompressible Cahn-Hilliard fluids and topological
  transitions}.
\newblock {\em {Proc. Royal Soc. A}}, {454}({1978}):{2617--2654}, {1998}.

\bibitem{PR05}
L.~Pareschi and G.~Russo.
\newblock Implicit-explicit {R}unge-{K}utta schemes and applications to
  hyperbolic systems with relaxation.
\newblock {\em J. Sci. Comput.}, 25(1/2):129--155, 2005.

\bibitem{ShenYang10}
Jie Shen and Xiaofeng Yang.
\newblock {Numerical Approximations of Allen-Cahn and Cahn-Hilliard Equations}.
\newblock {\em Discrete and Continuous Dynamical Systems}, 28(4):1669--1691,
  DEC 2010.

\bibitem{Shu09}
Chi-Wang Shu.
\newblock {High Order Weighted Essentially Nonoscillatory Schemes for
  Convection Dominated Problems}.
\newblock {\em SIAM Rev.}, 51(1):82--126, 2009.

\bibitem{Siano79}
Donald~B. Siano.
\newblock {Layered sedimentation in suspensions of monodisperse spherical
  colloidal particles}.
\newblock {\em J. Colloid and Interface Sci.}, 68(1):111--127, 1979.

\bibitem{Toro09}
E.~F. Toro.
\newblock {\em {Riemann solvers and numerical methods for fluid dynamics.}}
\newblock Springer, third edition edition, 2009.

\bibitem{Vollmayr-Lee-Rutenberg2003}
BP~Vollmayr-Lee and AD~Rutenberg.
\newblock {Fast and accurate coarsening simulation with an unconditionally
  stable time step}.
\newblock {\em {Phys. Rev. E}}, {68}({6, 2}), {2003}.

\bibitem{Yue04}
PT~Yue, JJ~Feng, C~Liu, and J~Shen.
\newblock {A diffuse-interface method for simulating two-phase flows of complex
  fluids}.
\newblock {\em Journal of Fluid Mechanics}, 515:293--317, SEP 25 2004.

\end{thebibliography}
\end{document}